\documentclass{article}
\usepackage[arrow,matrix,tips]{xy}
\usepackage{graphicx}
\usepackage{amsmath,amssymb}
\usepackage{amsfonts}
\usepackage{bm}
\def\t-hol{transversely holomorphic}

\def\L{{\mathcal L}}

\def\bgn{\begin}

\def\CL{\text{\rm CL}}

\def\E{\bold E}
\def\J{{\mathcal J}}
\def\L{{\mathcal L}}
\def\Ob{\text{\rm Ob}}
\def\1{{[1]}}
\def\2{{[2]}}
\def\3{{[3]}}
\def\({\left(}
\def\){\right)}
\def\s-circ{\,{\scriptstyle{\circ}}\,}
\def\<<{<\negthinspace \negthinspace<}
\def\Ad{\text{\rm Ad}}
\def\ad{\text{\rm ad}}
\def\noin{\noindent}

\def\bgn{\begin}
\def\cou{\ss\text{\rm co}}
\def\bE{\bold E}
\def\<{<\negthinspace \negthinspace <}

\title{Deformations of  generalized complex and \\generalized K\"ahler structures}

\author{Ryushi Goto 
}
\date{June 5 2007}

\begin{document}
\maketitle

\footnote{ 
$^{*}$Partly supported by the Grant-in-Aid for Scientific Research (C),
Japan Society for the Promotion of Science. 
}

\begin{abstract} In this paper
we obtain a stability theorem of generalized K\"ahler structures with one pure spinor under small
deformations of generalized complex structures.
(This is analogous to the stability theorem of K\"ahler manifolds by Kodaira-Spencer.)
We apply the stability theorem to a class of compact K\"ahler manifolds which admits deformations to
generalized complex manifolds and
obtain non-trivial
generalized K\"ahler structures on
Fano surfaces and toric K\"ahler manifolds.
In particular, we show that every nonzero holomorphic Poisson structure on a K\"ahler manifold induces deformations of nontrivial generalized K\"ahler structures.

 \end{abstract}

\tableofcontents

\pagenumbering{arabic}

\def\a{\alpha}
\def\b{\beta}
\def\bs{\backslash}
\def\e{\varepsilon}
\def\gam{\gamma}
\def\Gam{\Gamma}
\def\k{\kappa}
\def\del{\delta}
\def\lam{\lambda}
\def\ome{\omega}
\def\Ome{\Omega}
\def\sig{\sigma}
\def\Sig{\Sigma}
\def\eps{\varepsilon}
\def\thet{\theta}
\def\A{{\mathcal A} }
\def\E{{\mathcal E}}
\def\f{\bold {f}}
\def\G{G_2}
\def\D{\Bbb D}
\def\Diff{\text{\rm Diff}_0}
\def\Q{\Bbb Q}
\def\R{\Bbb R}
\def\C{\Bbb C}
\def\H{\text{\rm H}}
\def\Z{\Bbb Z}
\def\O{\Bbb O}
\def\P{\frak P}
\def\M{\frak M}
\def\N{\frak N}
\def\V{{\mathcal V}}
\def\X{{\mathcal X}}
\def\w{\wedge}
\def\({\left(}
\def\){\right)}
\def\G{G_2} 
\def\Cl{\frak Cl} 
\def\nab{\nabla}
\def\neg{\negthinspace}
\def\h{\hat}
\def\hrho{\hat{\rho}}
\def\wideh{\widehat}
\def\til{\tilde}
\def\wtil{\widetilde}
\def\ch{\check}
\def\l{\left}
\def\ol{\overline}
\def\pa{\partial}
\def\olpa{\ol{\partial}}
\def\r{\right}
\def\ran{\rangle} 
\def\lan{\langle}
\def\Spin{\text{Spin}(7)}
\def\ss{\scriptscriptstyle}
\def\trian{\triangle}
\def\hook{\hookrightarrow}
\def\hyper{hyperK\"ahler}
\def\Ker{\text{\rm Ker}}
\def\arrow{\longrightarrow}
\def\lrarrow{\Longleftrightarrow}
\def\CP{\dsize{\Bbb C} \bold P}
\def\SL{\ss{\text{\rm SL}}}
\def\Sym{\text{Sym}^3\,\h{\Ome}^1 _X}
\def\tilSym{\widetilde{ \text{Sym}^3}\,\h{\Ome}^1_X }
\def\SymS{\text{Sym}^3\,{\Ome}^1_S }
\def\W{\wedge_{\ome^0}}
\def\sw{{\ss{\w}}}
\def\noin{\noindent}
\def\bsh{\backslash}
\def\reg{\text{\rm reg }}
\def\:{\, :\,}
\def\CL{\text{\rm CL}}
\def\TT{T\oplus T^*}
\def\complex{generalized complex }
\def\K\"ahler{generalized K\"ahler}
\newtheorem{theorem}{Theorem}[section]
\newtheorem{definition}[theorem]{Definition}
\newtheorem{lemma}[theorem]{Lemma}
\newtheorem{proposition}[theorem]{Proposition}
\newtheorem{corollary}[theorem]{Corollary}
\newtheorem{example}[theorem]{Example}
\newtheorem{problem}[theorem]{Problem}
\newtheorem{conjecture}[theorem]{Conjecture}
\newenvironment{proof}{{\it proof}  }{q.e.d.}
\newenvironment{proof*}{ }{ q.e.d}


\large{
\setcounter{section}{-1}
\section{Introduction}

A notion of \complex structures was introduced  by Hitchin 
\cite {Hi1},
which interpolates between complex and symplectic structures.
An associated notion of \K\"ahler structures is developed by Gualtieri 
\cite{Gu1}.
Examples of \K\"ahler structures have been constructed by the reduction  \cite{B.G.C}, \cite {L.T}  
which is a generalization of the symplectic quotient construction.
Hitchin gave an explicit construction of
\K\"ahler structures on Del Pezzo surfaces 
by using holomorphic Poisson structures and 
suggested that 
\K\"ahler structures are related to holomorphic Poisson structures
\cite{Hi2}, \cite {Hi3}.

Kodaira and Spencer showed that
K\"ahler structures on compact complex manifolds are stable under sufficiently small deformations of complex structures \cite{K.S III}.
More precisely, if $V_0$ is a compact K\"ahler manifold, then any small deformation $V_t$ of $V_0$ is also a K\"ahler manifold.

 The purpose of this paper is to establish a stability theorem of \K\"ahler structures under small deformations of \complex structures. Applying the theorem, we shall obtain a systematic construction of non-trivial \K\"ahler structures which arise as deformations of
ordinary K\"ahler manifolds with holomorphic Poisson structures. 
The construction provides many examples by using both holomorphic Poisson structures and deformations of complex structures.
In our construction, it is intriguing to solve the problem of obstructions to deformations of \K\"ahler structures. 
Note that there exists an obstruction to deformations of \complex structures in general.
We assume that there exists a family of deformations of \complex structures on a \K\"ahler manifold $X$.
Then we apply the method in
\cite{Go2} and show that every obstruction to corresponding deformations of \K\"ahler structures vanishes. The method is a generalization of the one in unobstructed theorem of Calabi-Yau manifolds by Bogomolov-Tian-Todorov
\cite{Ti}, which is also applied to obtain unobstructed deformations and the local Torelli type theorem for
Riemannian manifolds with special holonomy group \cite{Go1}.
For the more precise statement of the stability theorem, we explain
\complex structures, \K\"ahler structures and in particular,
a relation to pure spinors.

The notion of \complex structures is based on an idea of
replacing the tangent bundle $T$ of a manifold with the direct sum of the tangent bundle $T$ and the cotangent bundle $T^*$. 
The fibre bundle of the direct sum $\TT$ admits an indefinite metric $\lan\,, \, \ran$ by which we obtain the fibre bundle 
SO$(\TT)$ with fibre the special orthogonal group. 
An almost \complex structure $\J$ is defined as a section of the fibre bundle SO$(\TT)$ with $\J^2=-$id, which gives rise to the decomposition
$(\TT)\otimes\C =L_\J\oplus \ol L_\J$, where $L_\J$ is $-\sqrt{-1}$-eigenspace of $\J$
and $\ol L_\J$ denotes its complex conjugate. 
Almost
\complex structures form an orbit of the action of the real Clifford group of the real Clifford algebra bundle $\CL$ with respect to $(\TT, \lan\,,\,\ran)$ ({\it cf.} \cite{Ch}).
A \complex structure is an almost \complex structure which is  integrable with respect to the Courant bracket.

A \K\"ahler structure is a pair $(\J_0,\J_1)$ consisting of commuting \complex structures
$\J_0$ and $\J_1$ which gives rise to a generalized metric $G:=-\J_0\J_1$.

The direct sum $\TT$ acts on differential forms on a manifold by the interior product and the exterior product.
For a differential form $\psi$, we define a subspace $L_\psi$
by $L_\psi:=\{\, E\in (\TT)\otimes\C\,|\, E\cdot \psi=0\, \}$.
A non-degenerate pure spinor is a differential form $\psi$ which gives a decomposition
$(\TT)\otimes\C=L_\psi\oplus \ol L_\psi$. Thus a non-degenerate pure spinor $\psi$
induces an almost \complex structure $\J_\psi$. It turns out that if a non-degenerate pure spinor $\psi$ is $d$-closed, then the induced structure $\J_\psi$ is integrable.
For a K\"ahler form $\ome$, the exponential
$e^{\sqrt{-1}\ome}$ is a non-degenerate pure spinor which induces the \complex structure $\J_\ome$.
From this point of view, we introduce
{\it a \K\"ahler structure with one pure spinor} as a pair $(\J, \psi)$ consisting of a
\complex structure $\J$ and a $d$-closed, non-degenerate pure spinor $\psi$ which induces
the \K\"ahler structure $(\J, \J_\psi)$.
Then we obtain the following stability theorem.\\ \\
{\bf Theorem 3.1 }{\it
Let $(\J, \,\psi)$ be a generalized K\"ahler structure with one pure spinor on a compact manifold $X$.
We assume that there exists an analytic family of \complex structures $\{\J_t\}_{t\in \trian}$ on $X$ with $\J_0=\J$ parametrized by the complex one dimensional open disk $\trian$ containing the origin $0$.
Then there exists an analytic family of \K\"ahler structures with one pure spinor
$\{\, (\J_t,\, \psi_t)\}_{t\in \trian'}$ with $\psi_0=\psi$ parametrized by a sufficiently small open disk $\trian'\subset \trian$ containing the origin.
}
\\ \\
An analytic family of \complex structures is a family of \complex structures $\{\J_t\}$ which depend analytically on
the parameter $t$ in $\trian$.
If the space of obstructions to deformations of \complex structures vanishes, 
then infinitesimal deformations generate an analytic family of deformations of 
\complex structures. 
It is remarkable that a holomorphic Poisson structure on a compact K\"ahler manifold gives the analytic family
of deformations of \complex structures which induces 
a family of deformations of non-trivial \K\"ahler structures.

In section 1, we present an exposition on \complex and \K\"ahler geometry. Preliminary results are collected in subsections 1-1 and 1-2 ({\it cf}. \cite{Gu1}, \cite{Gu2} and \cite{Hi1}).
In subsection 1-3, we introduce a \K\"ahler structure with one pure spinor and construct a differential complex $(K^\bullet,d)$ which is a subcomplex of the de Rham complex.
Applying the generalized Hodge decomposition \cite{Gu2},
we obtain an injective map from the cohomology
$H^*(K^\bullet)$ of the complex $(K^\bullet, d)$ to the de Rham cohomology group.
In section 2 we discuss deformations of \complex structures from the view point of pure spinors. The Maurer-Cartan equation naturally arises as the integrability of almost \complex structures. Further we show that an analytic family of \complex structures $\{\J_t\}_{t\in \trian}$ are described in terms of an analytic family of sections $a(t)$ of the real Clifford bundle $\CL^2$ with respect to $(\TT, \lan\,,\,\ran)$ which is the Lie algebra of
the Clifford group (conformal pin group). 
The exponential of sections $a(t)$ of $\CL^2$ is the family of sections of the Clifford group which acts on $\J_0$ by the adjoint action and we have
$$
\J_t=\Ad_{e^{a(t)}}\J_0.
$$%
We prove the stability theorem in section 3 in the sense of formal power series.
For the analytic family $a(t)$, we will construct a family of sections $b(t)$ of $\CL^2$ such that
\bgn{align}
&d\,(e^{a(t)}\, e^{b(t)}\, \psi_0)=0,\label{1}\\
&\Ad_{e^{b(t)}} \J_0=\J_0.\label{2}
\end{align}
It follows from the Campbell-Hausdorff formula \cite{Se} that we have a unique family $z(t) \in \CL^2$ with
$$
e^{z(t)}=e^{a(t)}\, e^{b(t)}.
$$%
Then from (\ref{1}), $e^{z(t)}\, \psi_0$ is a $d$-closed and non-degenerate pure spinor and we have $$\Ad_{e^{z(t)}}\J_0=\J_t,$$
from (\ref{2}). Since almost \K\"ahler structures also form the orbit of the action of the Clifford group,
it follows that $( \J_t,\,\,e^{z(t)}\, \psi)$ is a family of \K\"ahler structures with one pure spinor.
When we try to solve the equations (\ref{1}) and (\ref{2}), we encounter the class of obstruction
$[\wtil\Ob_k]\in H^2(K^\bullet)$ for each $k>0$. 
It turns out that each representative $\wtil\Ob_k$ is a $d$-exact differential form.
Since the cohomology group $H^2(K^\bullet)$ is embedded into the de Rham cohomology group, it follows that the class $[\wtil\Ob_k]$ vanishes and  
we obtain a solution $b(t)$ of the equations $(1)$ and $(2)$ as the formal power series.
Our solution $b(t)$ is not unique in general. 
A solution $b(t)$  together with $a(t)$ gives rise to a cohomology class of $H^1(K^\bullet)$ by the action on $\psi_0$.
We show that there exists a family of
solutions of the equations $(1)$ and $(2)$ which are locally parametrized by
the first cohomology group $H^1(K^\bullet)$ of the complex $(K^\bullet, d)$.
\\ \\
{\bf Theorem 3.2} {\it
Let $\{\J_t\}_{t\in \trian}$ and $\psi$ be as in theorem 3.1.
Then there is an open set $W$ in $H^1(K^\bullet)$ containing the origin such that there exists
a family of \K\"ahler structures with one pure spinor
$\{(\J_t, \,\psi_{t,s})\}$ with $\psi_{0,0}=\psi$
parametrized by $t\in \trian'$ and $s\in W$ in $H^1(K^\bullet)$.
Further if we denote by $[\psi_{t,s}]$ the de Rham cohomology class represented by $\psi_{t,s}$, then $[\psi_{t, s_1}]\neq [\psi_{t,s_2}]$ for $s_1\neq s_2$.
} \\ \\
In section 4, we will prove that the formal power series $b(t)$ converges and finish the proof of the stability theorem.
In section 5, we construct examples of \K\"ahler structures
on compact K\"ahler manifolds such as
Fano surfaces and toric manifolds.
Since there is no obstruction to deformations of \complex structures on any Fano surface, we can count the dimensions 
of deformations of \complex and \K\"ahler structures respectively.
We show that a holomorphic Poisson structure induces many interesting \K\"ahler structures.
If there is an action of a complex $2$-dimensional commutative Lie group which gives 
a nontrivial holomorphic Poisson structure on a compact K\"ahler manifold, 
then we obtain a family of deformations of nontrivial \K\"ahler structures.
It follows that every compact toric K\"ahler manifold admits nontrivial generalized K\"ahler structures.\\
There is a one to one correspondence between \K\"ahler structures and bihermitian structures \cite{Gu1}.
Then by using the stability theorem, it is shown that there exists a family of non-trivial bihermitian structures on 
every compact K\"ahler manifold $(X,\ome)$ with a non-zero holomorphic Poisson structure $\b$. 
Then we obtain an unobstructed deformations of complex structures whose infinitesimal deformation
is given by $\b\cdot\ome$ which is a $\ol\pa$-closed $T^{1,0}$-valued form of type $(0,1)$ given by the contraction of $\b$ by $\ome$. Thus we obtain\\ \smallskip
{\bf Theorem 3.2} \cite{Go3}
{\it Let $X$ be a compact K\"ahler manifold with a holomorphic Poisson structure $\b$. 
The class $[\b\cdot\ome]\in H^1(X, \Theta)$ gives rise to 
unobstructed deformations of complex structures.}
(see section 3 in \cite{Go3} for more detail).
\medskip

The author would like to thank Professor Fujiki and Professor
Namikawa for valuable discussions and suggestions. 
He wishes to thank Professor Hitchin for meaningful discussions.
After he posted his paper to Arxiv, he received a kind and sincere message from Professor Gualtieri.  
He is also grateful to Professor Yi Lin for his valuable message about the reduction.
}

\large{\numberwithin{equation}{section}
\section{Generalized complex and K\"ahler structures }
\subsection{\complex structures}
Let $\TT$ be the direct sum of the tangent bundle $TX$ and the cotangent bundle $T^*X$ on
a manifold $X$ of real $2n$ dimension.
Then there is a symmetric bilinear form $\lan\,,\,\ran$ on $\TT$ which is given by
\bgn{equation}\label{3}
\lan v+\theta, w+\eta\ran =\frac12\theta(w)+\frac12\eta(v),
\end{equation}
where $v, w\in TX$ and $\theta,\eta\in T^*X$.
Then we have the fibre bundle SO$(\TT)$ with fibre the special orthogonal group with respect to
$\lan\,,\,\ran$.
 We define {\it an almost generalized complex structure} $\J$ as a section of the bundle SO$(\TT)$ with $\J^2=-$id.
The direct sum $\TT$ acts on the differential forms $\w^\bullet T^*X$ by
the interior product and the exterior product,
\bgn{equation}\label{4}
(v+\theta)\cdot\a:= i_v\a+\theta\w\a,
\end{equation}
where $\a\in \w^\bullet T^*X$.
Let $\CL$ be the real Clifford algebra bundle of $\TT$ with respect to the bilinear form $\lan\,,\,\ran$.
Then from (\ref{3}) and (\ref{4}) we have the induced action of $\CL$ on differential forms $\w^\bullet T^*X$, which is the spin representation of $\CL$.
For a complex differential form $\phi$ we define a subspace  $L_\phi$ of $(\TT)\otimes\C$ by
\bgn{equation}
L_\phi:=\{\, E\in (\TT)\otimes\C\,|\,E\cdot\phi=0\,\}.
\end{equation}
A complex differential form $\phi$ is a (complex) {\it pure spinor} if $L_\phi$  is maximally isotropic, i.e., $2n$ dimensional.
A (complex) pure spinor $\phi$ is {\it non-degenerate} if we have the decomposition of
$(\TT)\otimes \C$ into $L_\phi$ and its complex conjugate $\ol L_\phi$,
\bgn{equation}
(\TT)\otimes\C=L_\phi\oplus \ol L_\phi.
\end{equation}
The decomposition (6) induces the almost generalized complex structure $\J_\phi$ which is defined by
\bgn{equation}
\J_\phi(E)=
\bgn{cases}
-\sqrt{-1} E, \quad &(E\in L_\phi),\\
\sqrt{-1}E,\quad &(E\in \ol L_\phi).
\end{cases}
\end{equation}
We call $\J_\phi$ the induced structure from the non-degenerate pure spinor $\phi$.

Let $\J$ be an almost \complex structure with the ${\scriptstyle -\sqrt{-1}}$-eigenspace $L_\J$.
Then we have the decomposition,
$(\TT)\otimes\C=L_\J\oplus \ol L_\J$.
We denote by $\CL^{[i]}$ the subbundle of $\CL$ of degree $i$. Then we identify the Lie algebra bundle {\it so}$\,(\TT)$ with $\CL^{[2]}$.
Under the identification {\it so}$\,(\TT)= \CL^{[2]}$, $\J$ acts on $\w^\bullet T^*X\otimes\C$ by
the spin representation. Then we have the eigenspace decomposition of $\w^\bullet T^*X\otimes\C$,
\bgn{equation}
\label{eq : 8}
\w^\bullet T^*X\otimes\C=U^{-n}\oplus U^{-n+1}\oplus \cdots \oplus U^{n-1}\oplus U^n,
\end{equation}
where $U^{k}$ denotes the eigenspace with eigenvalue $k\sqrt{-1}$.
The space $U^{-n}$ is a complex line bundle which we call the canonical line bundle of $\J$.
(We also denote it by $K_\J$).
Let $\w^k\ol L_\J$ be the $k$-th exterior product of $\ol L_\J$.
Then the eigenspace $U^{-n+k}$ is given by the action of $\w^k\ol L_\J$ on $K_\J$,
\bgn{equation}
U^{-n+k}=\w^k\ol L_\J \cdot K_\J.
\end{equation}
We denote by $\{(U_\a,\phi_\a)\}$ a trivialization of the line bundle $K_\J$, where 
$\{U_\a\}$ is a covering of $X$.  Each $\phi_\a$ is 
a non-vanishing section of $K_\J|_{U_\a}$ which is a non-degenerate pure spinor with the induced structure
$\J$.
Let $d$ be the exterior derivative and $E$ an element of $\CL^{[1]}\otimes\C=(\TT)\otimes\C$.
Then the anti-commutator $\{d, E\}:=dE+Ed$ acts on $\w^\bullet T^*X$.
 We have the derived bracket by the commutator of $\{d, E\}$ and $F$,
\bgn{equation}\label{eq: derived bracket}
[E,F]_d:=[\{d,E\}, F]
\end{equation}
 By skew-symmetrization of the derived bracket, we construct
the Courant bracket as 
\bgn{equation}
\label{eq : 10}
[E, F]_{\cou}:=\frac12[\{d,E\},F]-\frac12[\{d,F\},E].
\end{equation}
This is known as the derived bracket construction \cite{Kos}.
Note that if $E=v, F=w\in TX$, then the Courant bracket becomes the standard bracket of vector fields.
If the subbundle $L_\J$ is involutive with respect to the Courant bracket, then
$\J$ is {\it integrable}. A \complex structure is an almost \complex structure
which is integrable. 
The integrability of $\J$ is also given in terms of the corresponding pure spinor.
The following observation can be found in section 4.4 \cite{Gu1}.
\bgn{lemma}
Let $\phi$ be a non-degenerate pure spinor with the induced structure $\J_\phi$.
Then $\J_\phi$ is integrable if and only if there exists $E\in \CL^{[1]}\otimes\C=(\TT)\otimes\C$ such that
\bgn{equation}
d\phi+E\cdot\phi=0.
\end{equation}
\end{lemma}
To make the paper self-contained, we will give a proof.\\
\bgn{proof}
It suffices to show that $[E_1,E_2]_{\cou}\in L_\phi$ for $E_1,E_2\in L_\phi$.
It follows
\bgn{equation}
[\{d,E_1\},E_2]\phi=-E_2E_1d\phi.
\end{equation}
If we have $d\phi+E\cdot\phi=0$, then it follows
\bgn{align}
[\{d,E_1\},E_2]\phi=&E_2E_1E\phi,\\
=&\lan E_1,E\ran E_2\phi=0.
\end{align}
Hence from \text{\rm (\ref{eq : 10})},  we have $[E_1, E_2]_{\cou}\phi=0$. It implies that $L_\phi$ is involutive.
Conversely, assume that $\J$ is integrable.
From \text{\rm (\ref{eq : 8})},
$d\phi$ is decomposed into
\bgn{equation}
d\phi=\sum_{k=-n}^n(d\phi)^{[k]},
\end{equation}
where $(d\phi)^{[k]}\in U^{-n+k}$.
Then it follows that if $(d\phi)^{[k]}\neq 0$ for $k > -n+1$, then
there are $E_1, E_2$ such that
$[\{d,E_1\},E_2]\phi=-E_2 E_1d\phi\neq 0$.
Hence $d\phi\in U^{-n+1}$. It implies that
$(d\phi)=-E\cdot\phi$ for $E\in \CL^{[1]}\otimes\C$.
\end{proof}\\
If $\J$ is integrable, the image $d(U^k)$ is a subspace of the direct sum $U^{k-1}\oplus U^{k+1}$. Then $d$ is decomposed into
$\pa + \ol \pa$, 
$$
d\a =\pa \a+\ol \pa \a,
$$%
where $\pa\a\in U^{k-1}$ and $\ol\pa\a\in U^{k+1}$
for $\a\in U^k$.
There is a natural filtration of the even part of the real Clifford bundle $\CL$ ,
\bgn{align}
&\CL^0\subset \CL^2\subset\cdots.
\end{align}
We also have a filtration of the odd part of the real Clifford bundle,
\bgn{align}
&\CL^1\subset\CL^3\subset\cdots.
\end{align}
For instance, the first several ones are given by
\bgn{align*}
&\CL^0=C^\infty(X),\quad \CL^1=\CL^{[1]}=\TT, \\
&\CL^2=\CL^0\oplus\CL^{[2]}, \quad
\CL^3=\CL^{[1]}\oplus \CL^{[3]},
\end{align*}
where $\CL^{[i]}$ denotes the skew-symmetric subspace of $(\TT)$ in 
$\CL^i$.
The filtrations give rise to the filtration of bundles $\bE^k$ given by the action of $\CL^{k+1}$ on the canonical line bundle $K_\J$,
$$
\bE^k:=\CL^{k+1}\cdot K_\J,
$$%
where $\bE^k=\{0\}$ for $k< -1$. Note that $\bE^k$ is the complex vector bundle since 
$K_\J$ is the complex line bundle.
We change the degree of $\bE^\bullet$. For instance, $\bE^{-1}$ is the canonical line bundle $K_\J$ and
$\bE^0$ and $\bE^1$ are respectively written in the forms
\bgn{align}
\bE^0=&\{\, E\cdot \phi \,|\, E\in \CL^1, \phi \in K_\J\,\},\\
\bE^1=&\{\, a\cdot\phi\,|\, a\in \CL^2, \phi\in K_\J\,,\}.
\end{align}
Then $\bE^k$ is the direct sum in terms of $U^{-n+\bullet}$,
First four bundles are given by
\bgn{align}
&\bE^{-1}=U^{-n},\\
&\bE^{0}=U^{-n+1},\\
&\bE^1=U^{-n}\oplus U^{-n+2},\\
&\bE^2=U^{-n+1}\oplus U^{-n+3}.
\end{align}
Then $U^{-n+k}$ is given as the quotient bundle,
$$
U^{-n+k}=\bE^{k-1}/\bE^{k-3}.
$$%
It follows from $d=\pa+\ol\pa$ that $\bE^\bullet $ is invariant under the action of $d$.
Hence we have the differential complex $(\bE^\bullet, d)$,
\bgn{equation*}
\xymatrix{
0\ar[r]^d &\bE^{-1}\ar[r]^d &\bE^{0}\ar[r]^d&\bE^1\ar[r]^d
&\bE^2\ar[r]^d &\cdots.
 }
\end{equation*}
It is shown that the complex $(\bE^\bullet, d)$ is elliptic in
\cite{Go2}.
We denote by $H^k (\bE^\bullet)$ the $k$ th cohomology of the complex $(\bE^\bullet, d)$.

\subsection{generalized K\"ahler structures}
In this subsection, we use the same notation as in \cite{Gu2}.
Let $(\J_0, \J_1)$ be a pair of commuting \complex structures.
Then we define $\h G$ by the composition,
$$
\h G=-\J_0\J_1=-\J_1\J_0.
$$%
The symmetric bilinear form $G$ is given by
$G(E_1, E_2):= \lan \h GE_1, E_2\ran$ for
$E_1, E_2\in \TT$.
\bgn{definition}
A pair $(\J_0,\J_1)$ consisting of commuting \complex structures
is a \K\"ahler structure if the symmetric bilinear form $G$ is positive-definite.
\end{definition}
Let $U^p_{\J_i}$ be the eigenspace with respect to $\J_i$ for $i=0,1$.  Because we have the commuting pair $(\J_0, \J_1)$, we have the simultaneous
decomposition into eigenspaces,
$$
\w^\bullet T^*X\otimes\C= \oplus_{p,q} U^{p,q},
$$%
where $U^{p,q}=U^p_{\J_0}\cap U^q_{\J_1}$.
Then the image  of $ U^{p,q}$ by the exterior derivative $d$ is decomposed into four components
$U^{p+1,q+1}\oplus U^{p+1, q-1}\oplus U^{p-1,q-1}\oplus U^{p-1,q+1}$ which induces the decomposition of $d$,
$$
d=\ol\delta_++ \ol\delta_- +\delta_+ +\delta_-.
$$%
\[\xymatrix{
U^{p-1,q+1}\ar@{<-}[dr]_{\delta_-}&& U^{p+1,q+1}\ar@{<-}[dl]^{\ol\delta_+} \\
& U^{p,q}\ar[dr]_{\ol\delta_-}\ar[dl]^{\delta_+}&\\
U^{p-1,q-1}&&U^{p+1,q-1}
}
\]
\centerline{ Figure 1}
\subsection{\K\"ahler structures with one pure spinor}
We already see that a non-degenerate pure spinor $\psi$ is a differential form which induces
the almost \complex structure $\J_\psi$.
\bgn{definition}
Let $(\J, \psi)$ be a pair consisting of \complex structure $\J$ and a non-degenerate pure spinor $\psi$ with $d\psi=0$.
A pair $(\J, \psi)$ is a \K\"ahler structure with one pure spinor if the corresponding pair $(\J, \J_\psi)$ is a
\K\"ahler structure.
\end{definition}

We denote by $K^1$ the bundle $U^{0,-n+2}$ and define the
graded left module $K^\bullet$ generated by $K^1$ over the Clifford algebra $\CL$. We set $K^i=\{0\}$ for $i\leq 0$.
Then it follows
\bgn{align}
K^1=&U^{0,-n+2},\\
\label{eqn : 25}
K^2=&U^{1,-n+1}\oplus U^{-1,-n+1}\oplus U^{1,-n+3}\oplus U^{1,-n+3}.
\end{align}\\
\bigskip
\resizebox{1\hsize}{!}
{$ 
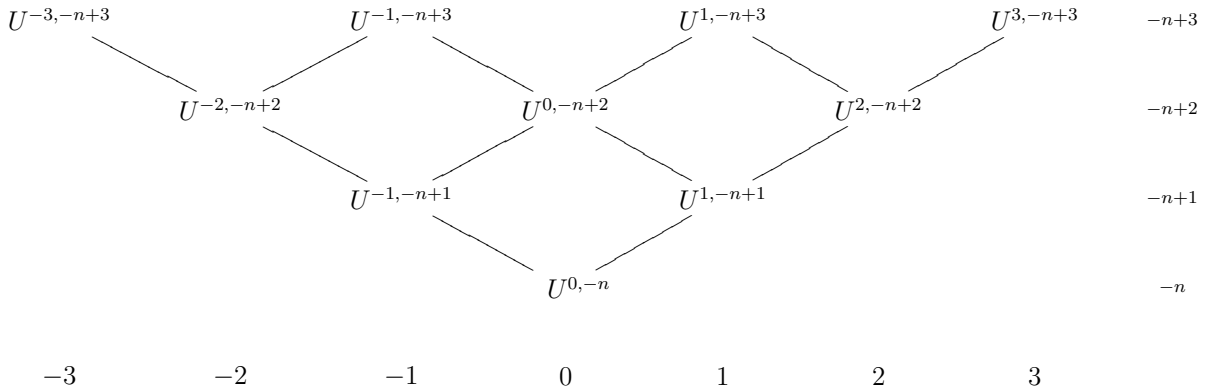
\xymatrix{
U^{-3,-n+3}\ar@{-}[dr]&&U^{-1,-n+3}\ar@{-}[dl]\ar@{-}[dr]&& U^{1,-n+3}\ar@{-}[dl]\ar@{-}[dr]&&U^{3,-n+3}\ar@{-}[dl]&{\scriptstyle -n+3}\\
&U^{-2,-n+2}\ar@{-}[dr]&& U^{0,-n+2}\ar@{-}[dl]\ar@{-}[dr]&&U^{2,-n+2}\ar@{-}[dl]&&\scriptstyle -n+2\\
&&U^{-1,-n+1} & &U^{1,-n+1}&&&\scriptstyle -n+1\\
&&&\quad U^{0,-n}\ar@{-}[ul]\ar@{-}[ur]&&&&\scriptstyle-n\\
-3&-2&-1&0 &1&2& 3&}
$}
\\\medskip
\large{
Then we have the following lemma from the decomposition of the exterior derivative $d$.
\bgn{lemma}
\label{lemma 2}
$(K^\bullet, d)$ is a differential complex.
\end{lemma}
Let $(\J, \psi)$ be a \K\"ahler structure with one pure spinor.
We denote by $a\cdot K_\J$ the action of $a\in \CL$ on the canonical line bundle $K_\J$.
We define a bundle $\ker^i$ by
\bgn{equation}\ker^i=\{\, a\in \CL^{i+1}\,|\, a\cdot K_\J=0\,\},
\end{equation}
for $i=0,1,2$.
We also define $\wtil\ker^i$ by using the filtration of $\CL$ and
$\bE^i:=\CL^{i+1}\cdot K_\J, $
\bgn{equation}
\wtil\ker^i=\{\, a\in \CL^{i+1}\,|\,a\cdot K_\J
\in \CL^{i-1}\cdot K_\J\,\}.
\end{equation}
Then we have
\bgn{lemma}
\bgn{equation}
U^{0,-n}\oplus U^{0,-n+2}=\{\, a\cdot\psi \, |\, a\in \wtil\ker^1\,\},
\end{equation}
\end{lemma}
\bgn{proof*}{\it proof of lemmas 1.5.}
The real bundle $\wtil\ker^1$ consists of linear combinations of the real part $E\cdot \ol{F}$ where $E\in L_{\J}$ and $\ol{F}\in \ol{L}_{\J}$. 
Since $E\cdot \ol{F}\psi \in U^{0,-n}\oplus U^{0,-n+2}$, it follows that 
$\wtil\ker^1\cdot \psi\in U^{0,-n}\oplus U^{0,-n+2}$.
Conversely it follows that $U^{0,-n}\oplus U^{0,-n+2}$ is generated by 
forms $(E\cdot \ol{F}+\ol{E}\cdot F)\psi$
and $\sqrt{-1}(E\cdot \ol{F}-\ol{E}\cdot F)\psi$
for $E\in L_{\J}$ and $\ol{F}\in \ol{L}_{\J}$. 
\end{proof*}
The bundle $K^2$ is also described in terms of $\ker^2$ and $\wtil{\ker}^2$,
\bgn{lemma}
\bgn{align}
K^2=&\{\, b\cdot\psi\,|\, b\in \ker^2\,\},\\
=&\{\, b\cdot\psi\, |\, b\in \wtil\ker^2\,\}.\nonumber
\end{align}
\end{lemma}
\bgn{proof*}{\it proof of lemma 1.6. }
We denote by $\wtil K^2$ the bundle $
\{\, b\cdot\psi\, |\, b\in \wtil\ker^2\,\}$.
Since $K^2$ is generated by $K^1$, we see that
\bgn{align}
\label{eqn : 30}
K^2\subset&\{\, b\cdot\psi\,|\, b\in \ker^2\,\}
\subset \wtil K^2.
\end{align}
The space $U^{3,-n+3}$ is given by $\w^3\ol L_\J\cdot\psi$.
Let $h$ be an element of $\w^3\ol L_\J$.
Then $h\cdot K_\J\in \CL^1\cdot  K_\J$ if and only if $h=0$. Since $\ker^2$ is real, $\wtil K^2$ does not contain the components $U^{3, -n+3}$ and $U^{-3,-n+3}$.
Hence it follows from (\ref{eqn : 25}) that $K^2=\wtil K^2$. We have the result from (\ref{eqn : 30}).
\end{proof*}
\bgn{lemma}
\label{lemma 3}
$(K^\bullet, d)$ is an elliptic complex for $i=1,2$.
\end{lemma}
\bgn{proof*}{\it proof of lemma \ref{lemma 3}.}
We will show that the symbol complex of the complex $(K^\bullet, d)$ is exact. It is sufficient to prove that if $u\w\a=0$ for non-zero
one form $u\in T^*$ and $\a\in K^i$ then $\a$ is given by
$\a=u\w\b$ for a $\b\in K^{i-1}$ for $i=1,2$.
We have the commuting generalized complex structures $\J$ and
$\J_\psi$ which act on $(\TT)\otimes\C$.
Then we have the simultaneous eigenspace decomposition,
\bgn{equation}
(\TT)\otimes\C = \ol L_+\oplus \ol L_-\oplus
L_+\oplus L_-,
\end{equation}
where $\ol L_+\oplus \ol L_-$ is $-{\scriptstyle\sqrt{-1}}$-eigenspace with respect to
$\J$ and $\ol L_+\oplus  L_-$ is $-{\scriptstyle \sqrt{-1}}$-eigenspace with respect to $\J_\psi$.
The non-zero element $u$ is decomposed into
\bgn{equation}\label{eqn : decomposition}
u=\ol u_+ +\ol u_- +u_+ + u_-,
\end{equation}
where $\ol u_\pm \in \ol L_\pm$ and $u_\pm\in L_\pm$.
Since $u\in T^*$,  we have $\lan u, u\ran=0$. Hence
\bgn{equation}\label{eq: 1.32}
0=\lan u, u\ran=\lan u_+,\ol u_+\ran+\lan u_-,\ol u_-\ran.
\end{equation}
The composition $\h G=-\J\J_\psi =-\J_\psi \J$ defines the generalized metric.
Since
$\h G(u_\pm+\ol u_\pm)=\pm (u_\pm+\ol u_\pm)$,
we have $(\pm1)\lan u_\pm, \ol u_\pm\ran > 0$.
In particular, it follows that
\bgn{equation}\lan u_\pm, \ol u_\pm \ran \neq 0,
\end{equation} 
because the generalized metric is positive-definite.
At first we consider the case  $K^1=U^{0,-n+2}$.
We assume that $u\w\a=0$ for non-zero $u\in T^*$ and
$\a\in U^{0,-n+2}$. Then it follows from the decomposition (\ref{eqn : decomposition}) that
\bgn{equation}
\ol u_\pm \cdot \a =0, \quad u_\pm\cdot \a=0.
\end{equation}
Then we have
\bgn{equation}
u_+\cdot\ol u_+\cdot\a =\lan u_+,\ol u_+\ran \a =0.
\end{equation}
Since $\lan u_+,\ol u_+\ran\neq 0$, we have $\a=0$.
In the case $K^2$, we assume that $u\w \a=0$ for non-zero $u\in T^*$ and $\a\in K^2$.
Form (\ref{eqn : 25}), we see that  $K^2\subset U_{\J_\psi}^{-n+1}\oplus U_{\J_\psi}^{-n+3}$.
Let $(\bE_\psi,d)$ be the differential complex defined by the action of $\CL$ on the canonical line bundle $K_{\J_\psi}$.
Since the complex $(\bold E_\psi,d)$ is elliptic, we have that
there exists $\til \b\in U_{\ss\J_\psi}^{\ss -n+2}$ such that
\bgn{equation}\label{eq: 1.36}
\a=u\w\til \b.
\end{equation}
We decompose $\til \b$ by
\bgn{equation}
\til \b=\til \b^{(2)}+\til \b^{(0)}+\til \b^{(-2)},
\end{equation}
where $\til \b^{(i)}\in U^{i,-n+2}$.
Then we define $\gam^{(\pm 1)}\in U^{\pm1, -n+1}$ by
\bgn{align}
&\gam^{(1)}=-\lan u_+, \ol u_+\ran^{-1} u_+\cdot\til\b^{(2)},\\
&\gam^{(-1)}=\lan u_-,\ol u_-\ran^{-1} \ol u_-\cdot\til\b^{(-2)}.
\end{align}
Then applying (\ref{eq: 1.32}) and (\ref{eq: 1.36}),
we obtain that
\bgn{align*}
u\w (u_-\cdot\gam^{(1)})=& (\ol u_+ + \ol u_-)\cdot u_-\cdot\gam^{(1)} \\
=&-(\ol u_+ + \ol u_-)\cdot u_-\cdot \lan u_+,\ol u_+\ran^{-1} u_+\cdot\til\b^{(2)}\\
=&(u_- + u_+)\cdot\til\b^{(2)}\\
=&u\w\til\b^{(2)},
\end{align*}
We also apply the similar method to  $\b^{(-2)}$,
then we have two equations
\bgn{align}
&u\w (u_-\cdot\gam^{(1)})=u\w\til\b^{(2)}\label{eqn : 41}
\\
\label{eqn : 42}
-&u\w(\ol u_+\cdot\gam^{(-1)})=u\w\til\b^{(-2)}.
\end{align}
We define $\b^{(0)}\in U^{0,-n+2}$ by
\bgn{equation}
\b^{(0)}=\til \b^{(0)}+u_-\cdot\gam^{(1)}-{\ol u}_+ \cdot\gam^{(-1)}.
\end{equation}
Then it follows from (\ref{eqn : 41}) and (\ref{eqn : 42}) that
\bgn{align}
u\w\b^{(0)}=&u\w\til \b^{(0)}+ u\w\til \b^{(2)}+u\w\til\b^{(-2)},\\
=&u\w\b=\a.
\end{align}
Hence the complex $(K^\bullet, d)$ is elliptic for $i=1,2$.
\end{proof*}

We denote by $H^i(K^\bullet)$ the $i$-th cohomology group of the complex $(K^\bullet,d)$.
The complex $(K^\bullet, d)$ is a subcomplex of the (full) de Rham complex $\{\cdots \overset d\to \w^\bullet T^*X\overset d\to \w^\bullet T^*X\overset d\to \cdots\} $. The cohomology group of the full de Rham complex is given by the full
de Rham cohomology group $H_{dR}(X):=\oplus_{i=0}^{2n} H^i(X,\C)$. Then we have the induced map
$p^i_{\ss K}\: H^i(K^\bullet)\to H_{dR}(X)$.
\bgn{lemma}
\label{lemma 4}
The map $p^i_{\ss K}\: H^i(K^\bullet)\to H_{dR}(X)$
is injective for $i=1,2$.
\end{lemma}
\bgn{proof*}{\it proof of lemma \ref{lemma 4}. }
Our proof is based on the generalized K\"ahler identities \cite{Gu2} (proposition 2),
\bgn{equation}\label{eqn : 46}
\ol\delta_+^*=-\delta_+, \quad \ol\delta_-^*=\delta_-,
\end{equation}
where the exterior derivative $d$ is given by
\bgn{equation}
\label{eqn : 47}
d=\ol\delta_+ +\ol\delta_-+\delta_+ +\delta_-,
\end{equation}
and $\ol\delta_\pm^*$ is the adjoint operator of $\ol\delta_\pm$ with respect to the generalized Hodge star operator.
Then the identities imply the equality of all available Laplacian,
\bgn{equation}\label{eqn : 48}
\trian_d=2\trian_{\ol\pa_\psi}=4\trian_{\ol\delta_\pm}
=4\trian_{\delta_\pm},
\end{equation}
where $\ol\pa_\psi=\ol\delta_+ +\delta_-$.
We obtain a $(p,q)$ decomposition for the de Rham cohomology of
any compact generalized K\"ahler manifold,
\bgn{equation}
H^\bullet(X,\C)=
\bigoplus_{\stackrel { |p+q|\leq n }{\ss p+q\equiv n \,(\text{\rm mod} 2)}} {\mathcal H}^{p,q},
\end{equation}
where ${\mathcal H}^{p,q}$ are $\trian_d$-harmonic forms in $U^{p,q}$.
At first we consider the cohomology $H^1(K^\bullet)$.
Let $\a$ be a $d$-closed element of $K^1$.
Then from (\ref{eqn : 47}) we have
\bgn{equation}
\ol\delta_\pm\a=0, \quad \delta_\pm \a=0.
\end{equation}
Then if follows from the generalized K\"ahler identities (\ref{eqn : 46}) that
\bgn{equation}
\ol\delta_+\a=0, \qquad\ol\delta_+^*\a=-\delta_+\a=0.
\end{equation}
Hence we have
\bgn{equation}
\trian_{\ol\delta_+}\a=(\ol\delta_+\ol\delta_+^*
+\ol\delta_+^*\ol\delta_+)\a=0.
\end{equation}
Then from (\ref{eqn : 48}),
$\a$ is $\trian_d$-harmonic and we have
\bgn{equation}
H^1(K^\bullet)\cong {\mathcal H}^{0,-n+2}.
\end{equation}
Hence we have the injection
$p^1_{\ss K}\: H^1(K^\bullet)\to H_{dR}(X)$.

In the case $H^2(K^\bullet)$, we use
the Green operators $G_{\ol\delta_\pm}$, 
$G_{\delta_\pm}$ and the Hodge decomposition of each $U^{p,q}$ by the elliptic operator $\trian_{\ol\delta_\pm}$.
We assume that $\a\in K^2$ is $d$-exact, i.e., $\a=d\b$.
Then it follows from $dd^{\J}$-lemma \cite{Gu2} that
we have an element of $\til \b\in U_{\J_\psi}^{-n+2}$ such that
\bgn{equation}
\a=d\til \b.
\end{equation}
(see the discussion \cite{Go2}.)
Then $\til \b$ is decomposed into the form,
\bgn{equation}
\til \b=\til \b^{(2)}+\til \b^{(0)}+\til \b^{(-2)},
\end{equation}
where $\til\b^{(i)}\in U^{i,-n+2}$.
We define $\gam^{(\pm1)}$ by
\bgn{align}
\gam^{(1)}=&\delta_+G_{\ol\delta_+}\til\b^{(2)}\\
\gam^{(-1)}=&\ol\delta_-G_{\delta_-}\til\b^{(-2)}
\end{align}
Then from the generalized K\"ahler identities (\ref{eqn : 46}) we have
\bgn{align}
&d\delta_-\gam^{(1)}=d\til\b^{(2)}
\label{eqn : 58}\\
-&d\ol\delta_+\gam^{(-1)}=d\til\b^{(-2)}
\label{eqn : 59}
\end{align}
We define $\b^{(0)}$ by
\bgn{equation}
\b^{(0)}=\til\b^{(0)}+\delta_-\gam^{(1)} -\ol\delta_+\gam^{(-1)}.
\end{equation}
Then it follows from (\ref{eqn : 58}) and (\ref{eqn : 59}) that
\bgn{align}
d\b^{(0)}=& d\til\b^{(0)} +d(\delta_-\gam^{(1)})-d(\ol\delta_+\gam^{(-1)})\\
=&d\til\b^{(0)}+d\til\b^{(2)}+d\til\b^{-2}\\
=&d\til\b=\a.
\end{align}
Hence every $d$-exact element $\a\in K^2$ is written as
\bgn{equation}
\a=d\b^{(0)},
\end{equation}
for $\b^{(0)}\in U^{0,-n+2}=K^1.$
It implies that the map $p^2_{\ss K}\: H^2(K^\bullet)\to H_{dR}(X)$
is injective.
\end{proof*}
}

\section{Deformations of generalized complex structures }
\large{
Let $\J$ be a \complex structure on a manifold $X$ with the maximally isotropic subspace $L(=L_\J)$ in $(\TT)\otimes\C$.
In the deformation theory of \complex structures developed in \cite{Gu1},
we will deform $L$ in the Grassmannian which consists of maximally isotropic subspaces. Then a small deformation of isotropic subspace is given by
\bgn{equation}\label{eq: 2.1}
L_\e:=(1+\e)L=\{ E+[E, \e]\,|\, E\in L\,\},
\end{equation}
for sufficiently small $\e\in \w^2\ol L$.
Then we have the decomposition $(\TT)\otimes \C$ into $L_\e$ and its complex conjugate $\ol L_\e$ which defines an almost \complex structure $\J_\e$ for $\e$.
The integrability of $\J_\e$ is equivalent to the one of almost Dirac structures in \cite{L.W.P}.
\bgn{theorem}\text{\rm (\cite{L.W.P})}
The structure $\J_\e$ is integrable if and only if $\e$ satisfies the generalized Maurer-Cartan equation,
\bgn{equation}
d_L\e+\frac12[\e,\e]_L=0,
\end{equation}
where $d_L\: \w^k\ol L\to \w^{k+1}\ol L$ denotes the exterior
derivative of the Lie algebroid  and
$[\,,\,]_L$ is the Lie algebroid bracket of $\ol L$, i.e.,
the Schouten bracket.
\end{theorem}
Let $\phi$ be a locally defined nowhere vanishing section of $K_\J$.
Then $\phi$ is a non-degenerate pure spinor which induces the structure $\J$.
The exponential $e^\e$ acts on $\phi$ and we have
the deformed non-degenerate pure spinor $e^\e\cdot\phi$
which induces $\J_\e$.
We already show that $\J_\e$ is integrable if and only if
the differential form $e^\e\phi$ satisfies
\bgn{equation}
de^\e\,\phi+ E_\e\cdot e^\e\,\phi=0,
\end{equation}
for $E_\e\in \CL^1\otimes\C$.
We will give another proof of theorem 2.1 from the view point of pure spinors. 
Our proof is suitable for our argument in this paper.\\
\bgn{proof}{\it of theorem 2.1.}
We recall the decomposition of differential forms,
\bgn{equation}
\w^\bullet T^*X\otimes\C=\bigoplus_{k=-n}^n U^k.
\end{equation}
Let $\pi_{U^{-n+3}}$ be the  projection to the component $U^{-n+3}$.
Since $\J_\e$ is integrable, we have
\bgn{equation}
de^\e\phi= -E_\e\cdot e^\e\phi,
\end{equation}
Let $\h E_\e$ be $e^{-\e}\,E\, e^\e\in \CL^1\otimes\C.$
Then by the left action of $e^{-\e}$ , we have
\bgn{equation}
e^{-\e}\,de^\e\,\phi= -\h E_\e\cdot\phi,
\end{equation}
We see that  $e^{-\e}d e^{\e}$ is a Clifford-Lie operator of order $3$ ({\it cf}. definition 2.2 in \cite{Go2}). It follows from definition that $e^{-\e}d e^{\e}$ is locally given by
the Clifford algebra valued Lie derivative,
$$
e^{-\e}d e^{\e}=\sum_i E_i \L_{v_i} +N_i,
$$%
where $\L_{v_i}$ is the Lie derivative by a vector filed $v_i$ and
$E_i\in \CL^1\otimes\C$, $N_i\in \CL^3\otimes\C$.
Thus $e^{-\e}de^\e\phi$ is an element of $U^{-n+1}\oplus U^{-n+3}$.
It implies that $\J_\e$ is integrable if and only if we have
$ \pi_{U^{-n+3}}\(e^{-\e}d e^\e\phi \)=0$.
The operator $e^{-\e}de^\e\phi$ is written in the form of power series ({\it cf}. lemma 2-7 in \cite{Go2})
\bgn{align}
\label{eqn : 1}
e^{-\e}de^{\e}\phi=&d\phi+[d,\e]\phi+\frac1{2!}[[d,\e],\e]\phi+\cdots,
\end{align}
We define $N(\e,\e)$ by
\bgn{equation}N(\e,\e):=[[d,\e],\e].\end{equation}
\bgn{lemma}
The operator $N(\e,\e)$ linearly acts on
$\w^\bullet T^*X$, which is  not a differential operator.
\end{lemma}
\bgn{proof}
We will show that
$[[d,\e_1],\e_2]f\a=f[[d,\e_1],\e_2]\a$
 for $\a\in\w^*T^*$ and a function $f$, where
 $\e_1, \e_2\in \w^2\ol L$.
It follows 
\bgn{align*}
[[d,\e_1],\e_2]f\a-&f[[d,\e_1],\e_2]\a\\
=(df)\e_1\e_2-&\e_1(df)\e_2-\e_2(df)\e_1+\e_2\e_1(df)\\
=(df)\e_1\e_2-&[\e_1,df]\e_2-[\e_2,df]\e_1+\e_2[\e_1,df]\\
-&(df)\e_1\e_2-(df)\e_2\e_1+[\e_2,(df)]\e_1\\
+&(df)\e_2\e_1   \\
=&[\e_2,[\e_1,(df)]].
\end{align*}
Since $\e_i\in\w^2\ol L$, we have $[\e_i,(df)]\in \ol L$. Hence
$$
[\e_i, [\e_j, (df)]]=0,
$$
for $i, j=1,2$.
Thus the result follows.
\end{proof}\\
The higher order terms of (\ref{eqn : 1}) are given by the adjoint action of $\e$ on $N(\e, \e)$ successively.
We define $\ad_\e^l N(\e,\e)$ by
$$
\ad_\e^l N(\e,\e):=[\ad_\e^{l-1} N(\e,\e), \e].
$$%
Hence we have
\bgn{align}
e^{-\e}d e^\e=&d\phi+[d,\e]\phi+\frac1{2!}N(\e,\e)\phi\\
+&\sum_{l=1}^\infty \frac1{(l+2)!}\,\ad_\e^l N(\e,\e).
\end{align}
Since $d_L$ is the exterior derivative of the Lie algebroid $\ol L$, we have the complex,
$$\xymatrix
{\cdots\ar[r]^{\ss{d_L}}&\w^p{\ol L}\ar[r]^{d_L}&\w^{p+1}{\ol L}\ar[r]^{\quad\ss d_L}&\cdots}.
$$
Then $d_L \e\in \w^3\ol L$ for $\e\in \w^2\ol L$ is given by
\bgn{lemma}
$$\pi_{U^{-n+3}}[d,\e]\phi=(d_L \e)\phi.$$
\end{lemma}
\bgn{proof}
Since we have $d\phi+E\phi=0$ for $E\in \ol L$, it follows
that
\bgn{equation}
\pi_{U^{-n+3}}(d+E)\e\phi=(d_L \e)\phi.
\end{equation}
Then we have
\bgn{align}
[d,\e]\phi=&d\e\phi-\e d\phi\\
=&d \e\phi +\e E\phi\\
=&d\e \phi+E \e\phi\\
=&(d+E)\e\phi.
\end{align}
Thus it follows
$$
\pi_{U^{-n+3}}[d,\e]\phi=(d_L\e)\phi.
$$
\end{proof}
\bgn{lemma}\label{lem: Schouten bracket}
The  Schouten bracket $[\e, \e]_L$ is given by
$$ [\e,\e]_L=N(\e,\e).$$
\end{lemma}
\bgn{proof}
Let $E_i$ be a section of $\TT$ for $i=1,2,3,4$.
In terms of the derived bracket $[E_i, E_j]_d=[\{ d, E_i\}, E_j]$ in (\ref{eq: derived bracket}),
the bracket $[[ d,\e_1],\e_2]$ is written as 
\bgn{align}
[[ d,\e_1],\e_2]=&-[E_1, E_3]_{d}E_2E_4 +[E_1, E_4]_{d}E_2E_3 \\
&+[E_2, E_3]_{d}E_1 E_4-[E_2, E_4]_{d} E_1E_3
\end{align}
for $\e_1=E_1E_2$ and $\e_2=E_3E_4$. 
Then the result follows. 
\end{proof}\\
Note that lemma \ref{lem: Schouten bracket} can be extended to higher order terms (see appendix).\\
We also have
\bgn{lemma}
$$\ad_\e^l N(\e,\e)=0,$$ for all $l\geq 1$.
\end{lemma}
\bgn{proof}
Since $N(\e,\e)\in \w^3\ol L$. It follows that
\bgn{equation}
[N(\e,\e), \e]=0.
\end{equation} Similarly we have
$\ad_\e^l N(\e,\e)=0$.
\end{proof}
 Then it follows from
lemma 2.3 and 2.4 that we have
\bgn{align}
\pi_{U^{-n+3}}\, e^{-\e}de^\e\phi= &d_L\e \phi+\frac1{2!}[\e,\e]_L\phi\\
=&\(d_L\e +\frac12[\e,\e]_L\)\phi.
\end{align}
Thus the equation
\bgn{equation}\pi_{U^{-n+3}}\, e^{-\e}de^\e\phi=0,\end{equation}
 is equivalent to  the Maurer-Cartan equation,
\bgn{equation}\(d_L\e +\frac12[\e,\e]_L\)=0.\end{equation}
Hence we have the result.
\end{proof}\\
Let $\e(t)$ be an analytic family of sections of $\w^2\ol L$. Then $\e(t)$ is written in the form of the power series in $t$,
\bgn{equation}
\e(t)=\e_1t+\e_2\frac{t^2}{2!}+\e_3\frac{t^3}{3!}+\cdots,
\end{equation}
where $t$ is a sufficiently small complex parameter.
Then $\e(t)$ gives deformations of almost \complex structures $\J_{\e(t)}$ by (\ref{eq: 2.1}).
The set of almost \complex structures forms an orbit of the adjoint action of SO$(\TT)$. 
The Lie algebra of SO$(\TT)$ is identified with $\w^2(\TT)$, which is the subspace $\CL^{[2]}$ of 
$\CL^2$. 
Thus $\J_{\e(t)}$ is written as $\J_{\e(t)}=\Ad_{e^{a(t)}}\J$ for $a(t)\in \w^2(\TT)$.
We denote by $(\w^2\ol L\oplus\w^2 L)^\R$ the real part of the bundle 
$(\w^2\ol L\oplus \w^2L)$ which is a subbundle of $\CL^2$.
Then we have 
\bgn{proposition}\label{prop: 2.6} 
There exists a unique analytic family $a(t)$ of sections of $(\w^2\ol L\oplus \w^2 L)^\R$ such that
\bgn{equation}
\J_{\e(t)}=\Ad_{e^{a(t)}}\J
\end{equation}
where we take sufficiently small $t$ if necessary.
\end{proposition}
\bgn{proof}
The action of $e^{\e(t)}$ on the canonical line bundle $K_\J$ defines a line bundle 
$e^{\e(t)}\cdot K_\J$. We also have a line bundle $e^{a(t)}\cdot K_\J$ by the action of $a(t)\in \CL^2$. 
The condition $ e^{\e(t)}\cdot K_\J=e^{a(t)}\cdot K_\J$ is equivalent to the condition $\J_{\e(t)}=\Ad_{e^{a(t)}}\J$.
Thus it suffices to construct a section $a(t)\in (\w^2\ol L\oplus \w^2 L)^\R$ which satisfies
\bgn{equation}\label{eq: 2-6-1}
(e^{-\e(t)}e^{a(t)})\phi\in K_J, \quad \text{\rm for all } \phi\in K_\J
\end{equation}
Given two differential forms $\a,\b$, if $\a-\b\in K_\J$, then we write it by 
$$
\a\equiv\b\, \,\,(\text{\rm mod} \,\,K_\J)
$$
Then the equation (\ref{eq: 2-6-1}) is written as 
$$
(e^{-\e(t)}e^{a(t)})\phi\equiv 0 \,\, \, (\text{\rm mod} \,\,K_\J)\quad \text{\rm for all } \phi\in K_\J
$$
We write $a(t)$ in the form of the power series in $t$,
\bgn{equation}
a(t)=a_1 t+a_2\frac{t^2}{2!}+\cdots,
\end{equation}
where $a_k$ is a section of $(\w^2\ol L\oplus \w^2 L)^\R$.
We denote by $(e^{-\e(t)}e^{a(t)})_{[k]}\phi$ the $k$ th term in $t$. Then the equation (\ref{eq: 2-6-1}) is reduced to infinitely many equations,
\bgn{equation}
(e^{-\e(t)}e^{a(t)})_{[k]}\phi\in K_\J,\quad \quad \text{\rm for all } \phi\in K_\J.
\end{equation}
We will show that there exists a solution $a(t)$ by induction on $k$.
For $k=1$, we have
\bgn{equation}
(e^{-\e(t)}e^{a(t)})_{[1]}\phi=-\e_1\phi+a_1\phi\in K_\J.
\end{equation}
Thus if we set
$a_1=\e_1+\ol\e_1$, then $(e^{-\e(t)}e^{a(t)})_{[1]}\phi=0\in K_\J$.
We assume that there are sections $a_1, \cdots, a_{k-1}\in (\w^2\ol L\oplus \w^2 L)^\R$
such that
\bgn{equation}\label{eq: 2-6-2}
\(e^{-\e (t)}e^{a(t)}\)_{[i]}\phi\in K_\J,
\end{equation}
for $\forall i<k$.
If follows from the
Campbel-Hausdorff formula that there exists $z(t)\in \CL^2\otimes\C$ such that
$
e^{-\e( t)}e^{a(t)}=e^{z(t)},
$
where
\bgn{equation}\label{eq: 2-6-3}
z(t)=-\e(t)+a(t)-[\e(t), a(t)]+\cdots.
\end{equation}
Thus our assumption (\ref{eq: 2-6-2}) is 
$
\(e^{z(t)}\)_{[i]}\cdot\phi \in  K_\J\,\,\text{\rm for all }i<k.
$
Since the degree of $z(t)$ is greater than and equal to $1$, we have $z(t)_{[1]}\cdot\phi\in  K_\J$ and 
it successively follows 
from our assumption that 
$z(t)_{[i]}\cdot\phi\in K_\J$, ($\forall i<k$).
Then we have 
\bgn{equation}
 (e^{z(t)} )_{[k]}\cdot\phi
 \equiv z(t)_{[k]}\phi\,\, (\text{\rm mod} \,\,K_\J)\quad \text{\rm for all } \phi\in K_\J
\end{equation}
Hence from (\ref{eq: 2-6-3}), there is a $H_k\in \CL^2\otimes\C$ such that
\bgn{align}
 (e^{z(t)})_{[k]}\cdot\phi
\equiv &\frac1{k!}a_k\phi -H_k\phi\,\,  (\text{\rm mod} \,\,K_\J)\quad \text{\rm for all } \phi\in K_\J
\end{align}
where $H_k$ is written in terms of  $a_1,\cdots , a_{k-1}$ and
$\e_1\cdots, \e_{k}$. Then there is a $\h H_k\in \w^2\ol L$ such that $\h H_k\phi -H_k\phi\in K_\J$.
Thus $a_k$ is defined as the real part of $(k!)\h H_k$ and we have
\bgn{align}
\frac1{k!}a_k\phi-H_k\phi\in K_\J.
\end{align}
Hence it follows
\bgn{equation}
 (e^{z(t)})_{[k]}\cdot\phi=\(e^{-\e(t)}e^{a(t)} \)_{[k]} \phi\in K_\J.
\end{equation}
Then we have a solution $a(t)$ as the formal power series.
It follows that the $a(t)$ is a convergent series which is a smooth section. 
Thus $a(t)$ is a unique section of $(\w^2\ol L\oplus\w^2 L)^\R$ with 
$\J_{\e(t)}=\Ad_{e^{a(t)}}\J$ which depends analytically on $t$.
\end{proof}
}

\large{
\section{Stability theorem of generalized K\"ahler structures }
We use the same notation as in sections 1 and 2.
\bgn{theorem}\label{stability}
Let $(\J, \,\psi)$ be a generalized K\"ahler structure with one pure spinor on a compact manifold $X$.
We assume that there exists an analytic family of \complex structures $\{\J_t\}_{t\in \trian}$ on $X$ with $\J_0=\J$ parametrized by the complex one dimensional open disk $\trian$ containing the origin $0$.
Then there exists an analytic family of \K\"ahler structures with one pure spinor
$\{\, (\J_t,\, \psi_t)\}_{t\in \trian'}$ with $\psi_0=\psi$ parametrized by a sufficiently small open disk $\trian'\subset \trian$ containing the origin.
\end{theorem}
Theorem \ref{stability} implies that \K\"ahler structures with one pure spinor are stable under
deformations of \complex structures. Theorem 3.1 is a generalization of
so called the stability theorem of K\"ahler structures due to Kodaira-Spencer.
We also obtain
\bgn{theorem}
Let $\{\J_t\}_{t\in \trian}$ and $\psi$ be as in theorem 3.1.
Then there is an open set $W$ in $H^1(K^\bullet)$ containing the origin such that there exists
a family of \K\"ahler structures with one pure spinor
$\{(\J_t, \,\psi_{t,s})\}$ with $\psi_{0,0}=\psi$
parametrized by $t\in \trian'$ and $s\in W$ in $H^1(K^\bullet)$.
Further if we denote by $[\psi_{t,s}]$ the de Rham cohomology class represented by $\psi_{t,s}$, then $[\psi_{t, s_1}]\neq [\psi_{t,s_2}]$ for $s_1\neq s_2$.
\end{theorem}
This section is devoted to prove theorem \ref{stability} and theorem 3.2.
Let $K_{\J_0}$ be the canonical line bundle with respect to $\J_0$.
We take a trivialization $\{U_\a, \phi_\a\}$ of $K_{\J_0}$, where
$\{U_\a\}$ is a covering of $X$ and $\phi_\a$ is a non-vanishing section
of $K_{\J_0}|_{U_\a}$ which induces the \complex structure $\J_0$.
Since $\J_0$ is integrable, we have $d\phi_\a +E_\a\phi_\a=0$ for $E_\a\in \CL^1\otimes\C|_{U_\a}$.
 It follows from section 2 that deformations $\{\J_t\}$ is given by an analytic family of global sections
 $a(t)\in \CL^2$ which is constructed from
 an analytic family of global sections $\e(t)\in \w^2\ol L$.
 Each section $a(t)$  gives  the non-degenerate pure spinor
 $ e^{a(t)}\phi_\a$ which induces the structure $\J_t$.
 Since $\J_t$ is integrable, we have
 \bgn{equation}
 de^{a(t)}\phi_\a+E_\a(t)e^{a(t)}\phi_\a=0.
 \end{equation}
 It follows from the left action of $e^{-a(t)}$
 \bgn{equation}
 e^{-a(t)}\,d\,e^{a(t)}\phi_\a+e^{-a(t)}E_\a(t)e^{a(t)}\phi_\a=0.
 \end{equation}
 We define $\til E_\a(t)$ by
 \bgn{equation}
 \wtil E_\a(t)=e^{-a(t)}E_\a(t)e^{a(t)}\in (\TT)|_{U_\a}=(\CL^1)|_{U_\a}.
 \end{equation}
 Then we have
 \bgn{equation}\label{eq : 100}
  e^{-a(t)}\,d\,e^{a(t)}\phi_\a+\wtil E_\a(t)\phi_\a=0
 \end{equation}
 Hence it follows
 \bgn{equation}
 (e^{-a(t)}\,d\,e^{a(t)})\phi_\a\in \bE^0_{\J_0}|_{U_\a}=\{\, E\cdot\phi_\a\,|
 \, E\in \CL^1|_{U_\a}\,\}.
 \end{equation}
 Since $e^{-a(t)}\,d\,e^{a(t)}$ is a Clifford-Lie operator of order $3$ ({\it cf}. definition 2.2 in \cite{Go2}),
 it follows that $e^{-a(t)}\,d\,e^{a(t)}$ is locally written in terms of the Lie derivative and the Clifford algebra,
 \bgn{equation}\label{eq: Clifford-Lie operator}
 e^{-a(t)}\,d\,e^{a(t)}=\sum_i E_i\L_{v_i}+N_i,
 \end{equation}
 where $E_i\in \CL^1$, $v_i\in T$ and $N\in \CL^3$.
 Then we have 
 \bgn{lemma} There is a section $a_i\in \CL^2$ such that 
 \bgn{align}
& \L_{v_i}\phi_\a \equiv a_i\cdot \phi_a\, \,
\text{\rm mod}\,(K_{\J_0}) \label{eq: mod K}\\
 &\L_{v_i}\psi=a_i\cdot\psi,
 \end{align}
 for each vector field $v_i$,
 where the equation (\ref {eq: mod K}) implies that 
 $$\L_{v_i}\phi_\a - a_i\cdot \phi_a= \rho_\a\phi_\a$$ for a function $\rho_a$.
 \end{lemma}
 \bgn{proof} The set of almost \K\"ahler structures with one pure spinor forms an orbit 
 under the diagonal action of the Clifford group whose Lie algebra is given by CL$^2$. 
 Thus small deformations of the structures are given by the exponential 
 action of CL$^2$. Let $f_t$ be the one parameter subgroup of diffeomorphisms 
 defined by the vector field $v$, i.e., 
 $$
 \frac{d}{dt} f_t |_{t=0}= v.
 $$
 Since the set of almost \K\"ahler structures with one pure spinor is invariant under the action of diffeomorphisms, there is a section $a(t)\in \CL^2$ with $a(0)=0$ such that 
 $$
 (f_t^*\J_0,\, f_t^*\psi) = (\Ad_{e^{a(t)}}\J_0, \,\,e^{a(t)}\cdot\psi).
 $$
 By differentiating with respect to $t$, we have 
 $$
 (\L_v\J_0 , \L_v\psi) = ( [ a, \J_0], \, a\cdot\psi), 
 $$
 where $a=\frac{d}{dt} a(t)|_{t=0}$.
 Since $f_t^*\phi_a$ and $e^{a(t)}\phi_\a$ induce the same \complex structure $\Ad_{e^{a(t)}}\J_0$, we have 
 $$
 f_t^*\phi_\a=e^{\rho(t)} e^{a(t)}\phi_\a,
 $$
 for a function $\rho(t)$ with $\rho(0)=0$.
 Then we have 
 \bgn{align}
& \L_{v}\phi_\a \equiv a\cdot \phi_a\, \,
\text{\rm mod}\,(K_{\J_0}) \\
 &\L_{v}\psi=a\cdot\psi.
 \end{align}
 \end{proof}\\
 Hence 
 it follows from (\ref{eq: Clifford-Lie operator})
 that there exists a section $h_\a\in \CL^3|_{U_\a}$ such that
 \bgn{align}
& (e^{-a(t)}\,d\,e^{a(t)})\phi_\a \equiv h_\a\cdot\phi_\a  \,\,
\text{\rm mod}\,(\CL^1\cdot K_{\J_0})\label{104}\\
& (e^{-a(t)}\,d\,e^{a(t)})\psi=h_\a\cdot\psi.
 \end{align}
 Let $K^\bullet$ be the graded left module generated by $U^{0,-n+2}$ over the Clifford algebra $\CL$.
 as in section 1.3.
 The exterior derivative $d$ gives rise to the differential complex :
 \bgn{equation}
 0\to K^1\to K^2\to \cdots.
 \end{equation}
 Then we see that
 $K^2$ is given by
 \bgn{equation}
 K^2=U^{1,-n+1}\oplus U^{-1,-n+1}\oplus U^{1,-n+3}\oplus U^{-1,-n+3}.
 \end{equation}
 We define a vector bundle $\ker^i$ by
 \bgn{equation}
 \ker^i=\{\, a\in \CL^{i+1}\,|\, a\cdot \phi_\a=0\,\},
 \end{equation}
 for $i=1,2$.
In section 1, 
 we define a bundle $\wtil \ker^i$ by
 \bgn{equation}
 \wtil\ker^i=\{\, a\in \CL^{i+1}\,|\, a\cdot\phi_\a\in
 \CL^{i-1}\cdot K_{\J_0}\,\}.
 \end{equation}
 The $\wtil\ker^i$ gives the bundle
 \bgn{equation}
 \til K^i=\{\, a\cdot\psi\, |\, a\in \wtil\ker^i\,\}.
 \end{equation}
 In section 1.3 we also have
 \bgn{align}
& \til K^1=U^{0,-n}\oplus U^{0,-n+2},\\
&\til K^2=K^2.
 \end{align}
 Hence $K^1$ is the subbundle of $\til K^1$,
 \bgn{equation}
 K^1\subset \til K^1.
 \end{equation}
 \bgn{proposition}\label{3.3}
 \[e^{-a(t)}\,d\,e^{a(t)}\psi\in K^2.\]
 \end{proposition}
 \bgn{proof}
 It follows from (\ref{104}) that there exists $h_\a\in \CL^3|_{U_\a}$ for each $\a$ such that
 \bgn{align}
 & e^{-a(t)}\,d\,e^{a(t)}\phi_\a\equiv h_\a\cdot\phi_\a \label{eq: 3.11}\,\,\,
 \text{\rm mod}\,(\CL^1\cdot K_{\J_0})\\
 & e^{-a(t)}\,d\,e^{a(t)}\psi=h_\a\cdot\psi,
 \end{align}
 where (\ref{eq: 3.11}) implies that 
 there is a section $F_\a\in \TT$ such that 
 $e^{-a(t)}\,d\,e^{a(t)}\phi_\a - h_\a\cdot\phi_\a= F_\a\cdot\phi_\a$.
 Since $\J_t$ is integrable, from (\ref{eq : 100}) we have
 \bgn{equation}e^{-a(t)}\,d\,e^{a(t)}\phi_\a=
  -\wtil E_\a(t)\cdot\phi_\a
 \in \CL^1\cdot K_{\J_0}|_{U_\a}.
 \end{equation}
 Hence it follows $h_\a\in \wtil\ker^2$ and we have
 \bgn{equation}
  e^{-a(t)}\,d\,e^{a(t)}\psi=h_\a\cdot\psi \in \wtil K^2=K^2.
 \end{equation}
\end{proof}\\ \\
\bgn{proof*} 
{\it proof of theorem 3.1 and theorem 3.2}
 We will construct a smooth family  $b(t)$ of sections of $\ker^1$ such that
  \bgn{equation}
  d(\,e^{a(t)} \,e^{b(t)}\,\psi)=0.
  \end{equation}
  Then it follows from the Campbel-Haudorff formula that there exists $z(t)\in \CL^2$ such that
\bgn{equation}
e^{z(t)}=e^{a(t)}e^{b(t)}.
\end{equation}
Explicitly, the first five components of $z(t)$ are given by
\bgn{align}
z(t)=&a(t)+b(t)+\frac12 [a(t), b(t)]\\
&+\frac1{12}[x,[x,y]]+\frac1{12}[y,[y,x]]+\cdots,
\end{align}
({\it cf}. \cite{Se}.)
  Since $b(t)\in \ker^1$, we have
  \bgn{align}
  e^{z(t)}\phi_\a=&e^{a(t)}e^{b(t)}\phi_\a\\
  =&e^{a(t)}\phi_\a.
  \end{align}
 It implies that $e^{z(t)}\phi_\a$ induces the same deformations
 $\J_t$ as before and the pair $(\J_t,\, e^{z(t)}\psi)$
 gives deformations of generalized K\"ahler structure with
 one pure spinor.
 Consequently the equation we must solve is that
 \bgn{align*}
   \qquad d\,e^{a(t)} \,e^{b(t)}\,\psi=0, \quad b(t)\in \ker^1.&
  \tag{\rm eq}
 \end{align*}
 The section $a(t)$ is written as the power series,
 \bgn{equation}
 a(t) =a_1 t+a_2\frac{t^2}{2!}+a_3\frac{t^3}{3!}+\cdots,
 \end{equation}
 where $a_i\in \CL^2$.
 We shall construct a solution $b(t)$ as the formal power series,
 \bgn{equation}
 b(t)=b_1t +b_2\frac{t^2}{2!}+b_3\frac{t^3}{3!}+\cdots,
 \end{equation}
 where $b_i\in \ker^1$.
 The $i$-th homogeneous part of the equation (eq) in $t$ is denoted by
 \bgn{equation*}
   \(d\,e^{a(t)} \,e^{b(t)}\,\psi\)_{[i]}=0, \quad b(t)\in \ker^1.
  \tag{eq$_{[i]}$}
 \end{equation*}
 Thus in order to obtain a solution $b(t)$, it suffices to determine
 $b_1,\cdots, b_i$ satisfying (eq)$_{[i]}$ by induction on $i$.
 In the case $i=1$, we have
 \bgn{align}
 \big(\,d\,e^{a(t)}\,e^{b(t)}\,\big)_{[1]}\psi=&da_1\psi +db_1\psi,\\
 =&[d,a_1]\psi+db_1\psi=0.
 \end{align}
 From proposition \ref{3.3} we have $\( e^{-a(t)}\,d\,e^{a(t)}\psi\)_{[1]}= [d,a_1]\psi\in K^2$.
 Since $da_1\psi=[d, a_1]\psi\in K^2$ is a $d$-exact differential form, $d a_1\psi$
 defines a class of cohomology $[\wtil\Ob_1]$ in $H^2(K^*)$ whose image vanishes in the de Rham cohomology group $H_{dR}(X)$.
 Since the map $p^2_K\: H^2(K^\bullet)\to H_{dR}(X)$ is injective, 
 it follows that $[\wtil\Ob_1]=0$. Thus we have a solution $b_1\in\ker^1$ which is given by
 \bgn{equation}
 b_1\psi=-d^*G_K(da_1\psi)\in K^1,
 \end{equation}
 where $d^*$ is the adjoint operator and $G_K$ is the Green operator of the complex $(K^*,d)$ with respect to a metric.
 Further for each representative $s$ of the first cohomology group $H^1(K^\bullet)$, we have a solution $b_{1, s}$ which is defined by
 \bgn{equation}
 b_{1,s}\psi=-d^*G_K(da_1\psi)+s.
 \label{eq : 131}
 \end{equation}

Assume that we already have $b_1,\cdots ,b_{k-1}\in \ker^1$ such that
  \bgn{equation}\label{eq : 132}
  \(de^{a(t)}e^{b(t)}\psi \)_{[i]}=0,
  \end{equation}
  for all $i<k$.
  From the Campbel-Hausdorff formula we have
  \bgn{equation}
  e^{z(t)}=e^{a(t)}e^{b(t)}.
  \end{equation}
  Hence it follows from our assumption (\ref{eq : 132})
  \bgn{align}
  \big(e^{-z(t)}\,d\,e^{z(t)}\big)_{[k]}\psi=&\sum_{\stackrel {i+j=k}{\scriptscriptstyle i,j\geq 0}}\(e^{-z(t)}\)_{[i]}\( de^{z(t)}\)_{[j]}\psi
  \nonumber\\
  \label{eq : 134}\\
  =&\(de^{z(t)}\)_{[k]}\psi. \nonumber
  \end{align}
Since $(e^{-z(t)}de^{z(t)})$ is given by
\bgn{equation}
(e^{-z(t)}de^{z(t)})=d+ [d, z(t)]+\frac1{2!}\left[ [d,z(t)], z(t)\right]+\cdots,
\end{equation}
the left hand side of (\ref{eq : 134}) is written as
\bgn{equation}
 (e^{-z(t)}de^{z(t)})_{[k]}\psi=\frac1{k!}d b_k \psi+\frac1{k!}da_k\psi
 +\Ob_k, \nonumber
\end{equation}
where $\Ob_k$ is the higher order term which is determined by
$a_1,\cdots , a_{k-1},$ and $ b_1,\cdots b_{k-1}$.
We define $\wtil{\Ob_k}$ by
\bgn{equation}
\wtil{\Ob_k}=\frac1{k!}da_k\,\psi+\Ob_k.
\end{equation}
Then the (eq)$_{[k]}$ is reduced to
\bgn{equation}
\frac1{k!} db_k \psi
 +\wtil{\Ob_k}=0, \qquad (b_k\in \ker^1)\nonumber
\end{equation}
From (\ref{eq : 100}), we have
\bgn{align}
e^{-z(t)}\, d\, e^{z(t)}\phi_\a=&e^{-b(t)}\, e^{-a(t)}\, d\,e^{a(t)}\, e^{b(t)}\phi_\a \nonumber\\
=&-\(e^{-b(t)}\til E_\a(t) e^{b(t)}\)\phi_\a\in \CL^1\cdot K_{\J_0}
\nonumber
\end{align}
Thus it follows from the same argument as in proposition 3.3 that we have
  \bgn{equation}(e^{-z(t)}de^{z(t)})\psi \in K^2.
  \end{equation}
It follows from (\ref{eq : 134}) that $\wtil{\Ob_k}\in K^2$ is $d$-exact. It implies that $\wtil{\Ob_k}$
gives rise to the class of the cohomology $[\wtil{\Ob_k}]\in H^2(K^*)$
with $p^2_K(([\wtil{\Ob_k}])=0$.
Since $p^2_K$ is injective from lemma 1.8, we have $[\wtil{\Ob_k}]=0$.
Thus $b_k\in \ker^1$ is given by
\bgn{equation}
\frac1{k!}b_k\psi=-d^*G_K(\wtil{\Ob_k})\in K^1,
\end{equation}
where $d^*$ is the adjoint operator and $G_K$ is the Green operator of the complex $(K^\bullet, d)$.
Hence it follows from the induction that we have the solution $b(t)$ of the equation (eq) as the formal power series.
As we see (\ref{eq : 131}), we obtain the family of sections $b_{1,s}$ parametrized by $s\in H^1(K^\bullet)$ which gives rise to a family $b(t, s)$ of solutions.
A family of non-degenerate pure spinor $\{\psi_{t,s}\}$ are
constructed as $e^{b(t,s)}\cdot\psi_0$.
Since the map $p^1_K\: H^1(K^\bullet)\to H_{dR}(X)$ is injective,
we have $[\psi_{t,s_1}]\neq [\psi_{t, s_2}]\in H_{dR}(X)$ for $s_1 \neq s_2$.
In section 4 we show that the formal power series $b(t)$ converges.
\end{proof*}
}

\large{
\section{The convergence }
This section is devoted to show that both power series $b(t)$ and $z(t)$ in section 3
are convergent series. We will use a similar method as in \cite{Ko} which apply the elliptic estimate 
of the Green operator. 
However we must develop an estimate of the obstruction $\Ob$ in section 3 which includes the higher order term.
We will use the induction on the degree $k$.
At first we will estimate the first terms $b_1$ and $z_1$ of power series $b(t)$ and $z(t)$.
We assume that $b(t)$ and $z(t)$ satisfy the inequality (\ref{155}) and (\ref{156}) respectively.
Then we will show that $b(t)$ satisfies the inequality (\ref{145}) and then obtain
the inequality (\ref{146}).\\
We shall fix our notation. We denote by $\|f\|_s=\|f\|_{C^{s,\a}}$ the H\"older norm of a section $f$ of a bundle with respect to a metric.
Then we have an inequality,
$$
\|fg\|_{s}\leq C_s\|f\|_{s}\, \|g\|_{s},
$$
where $f, g$  are sections and $C_s$ is a constant.
We have the elliptic complex $(K^\bullet, d)$ in section 1 and
we use the Schauder estimates of the elliptic operators with respect to the complex $(K^\bullet, d)$
with a constant $C_K$.
Let $P(t)$ be a formal power series in $t$.
We denote by $(P(t))_{[k]}$ the $k$ th coefficient of $P(t)$ and
Given two power series $P(t)$ and $Q(t)$,
if $ (P(t))_{[k]} < (Q(t))_{[k]} $ for all $ k$, we denote it by
$$P(t) \<Q(t).$$
For a positive integer $k$,
if $(P(t))_{[i]}<(Q(t))_{[i]}$ for all $i\leq k$, then we write it by
$$P(t)\underset k{\<}Q(t).$$
We also consider a formal power series  $f(t)$ in $t$ whose coefficients are sections of a bundle.
Then we put $\|f(t)\|_s =\sum_i \|\(f(t)\)_{[i]}\|_s t^i$.
We define a convergent power series $M(t)$ by
$$
M(t)=\sum_{\nu=1}^\infty\frac{1}{16c}\frac{(ct)^\nu}{\nu^2}=\sum_{\nu=1}^\infty M_\nu t^\nu.
$$
In \cite{Ko}, it turns out that the series $M(t)$ satisfies
\bgn{lemma}\label{lem: 4.1}
$$
M(t)^2 \<\frac 1c M(t).
$$
\end{lemma}
We put $\lam=\frac 1c$.  Then it follows from lemma 4.1 that
$$
\frac 1{l!}M(t)^l\<\frac1{l!}\lam^{l-1}M(t)=
\frac{\lam^l}{l!}\frac 1{\lam}
M(t).
$$
Hence we have
\bgn{lemma}\label{lem: 4.2}
$$ e^{M(t)}\<\frac 1\lam e^\lam M(t).
$$
\end{lemma}
As in section 3, the power series $z(t)$ is defined by the Campbel-Hausdorff formula,
$$
e^{z(t)}=e^{a(t)}e^{b(t)},
$$
where
\bgn{align}
z(t)=&\sum_{l=0}^\infty\frac {t^l}{l!} z_k,\\
e^{z(t)}=&\sum_{j=0}^\infty\frac {1}{j!}z(t)^j \\
=&1+z(t)+\frac1{2!}z(t)^2+\cdots.\nonumber
\end{align}
The power series $a(t)$ is the convergent series which induces deformations of
\complex structures $\{\J_t\}$ defined in proposition \ref{prop: 2.6}.
The norm of $a(t)$ is written as
$$
\|a(t)\|_s=\sum_{l=1}^\infty\frac1{l!}\|a_k\|_st^l.
$$%
Then we can assume that $\|a(t)\|_s$ satisfies
\bgn{equation}\label{142}
\|a(t)\|_s\<K_1M(t),
\end{equation}
for a non-zero constant $K_1$ and $\lam$ if we take $a(t)$ sufficiently small.
We will show that there exist constants $K_1$, $K_2$ and $\lam$ such that
we have the following inequalities,

\bgn{align}
&\|b(t)\|_s{\<}K_2M(t),\label{143}\\
&\|z(t)\|_s{\<}M(t) \label{144}
\end{align}
for sufficiently small $a(t)$.
Note that $K_1$, $K_2$ and $\lam$ are determined by
$a(t)$, $\J$ and $\psi$ which do not depend on
$b(t)$ and $z(t)$.
The inequalities (\ref{143}) and (\ref{144}) are reduced to the infinitely many inequalities on degree $k$
\bgn{align}
&\|b(t)\|_s\underset{k}{\<}K_2M(t),\label{145}\\
&\|z(t)\|_s\underset{k}{\<}M(t) \label{146}
\end{align}
We will show both inequalities (\ref{145}) and (\ref{146})
by the induction on $k$. In this section we denote by $C_i$  constants which do not depend on $z(t)$, $b(t)$ and $k$ but depend on $a(t)$, $\J$ and $\psi$.
For $k=1$,
as in section 3,
$b_1\psi$ satisfies the equation,
$$
d b_1\psi +da_1\psi=0,\qquad (b_1\psi\in K^1)
$$
Then $b_1\psi$ is given by
\bgn{equation}
b_1\psi=-{d}^*{G_K}(d a_1\psi),
\end{equation}
where $d^*$ is the adjoint operator and $G_K$ is the Green operator of the complex $(K^\bullet, d)$.
If follows from the Schauder estimate of the elliptic operators that
\bgn{align} \label{148}
\|b_1\psi\|_s&\leq C_{K}\|a_1\psi\|_s\leq C_{K}C_s\|a_1\|_s\|\psi\|_s\\
&\leq \frac 1{16}C_1K_1,\nonumber
\end{align}
where $\|a_1\|_s\leq K_1 M_1=\frac {K_1}{16}$ and $C_1=C_KC_s\|\psi\|_s$.

We can define $b_1$ as a section of the real part of $\ol L_+ L_-$.
Then we have
\bgn{equation}\label{149}
\|b_1\|_s\leq C_2\|b_1\psi\|.
\end{equation}
Substituting (\ref{148}) into (\ref{149}), we have
\bgn{equation}\label{150}
\|b_1\|_s\leq \frac 1{16}C_1C_2K_1=M_1C_1C_2K_1
\end{equation}
Thus if we take $K_2$ with
$C_1C_2K_1< K_2$, then we have
\bgn{equation}\|b_1\|_s\leq K_2 M_1,\label{151}
\end{equation}
Since $z_1=a_1+b_1$, if we take $K_1$ and $K_2$ satisfying $K_1+K_2<1$, we have
\bgn{align}
\|z_1\|_s &\leq \|a_1\|_s+\|b_1\|_s\\
&\leq M_1K_1+M_1K_2\\
&=(K_1+K_2)M_1<M_1  \label{154}
\end{align}
It follows from (\ref{151}), (\ref{154}) that we have inequalities (\ref{145}) and (\ref{146}) for
$k=1$.
\\
We assume that the following inequalities hold
\bgn{align}
&\|b(t)\|\underset{k-1}{\<}K_2M(t) \label{155}\\
&\|z(t)\|\underset{k-1}{\<}M(t)  \label{156}.
\end{align}
Let $\Ob_k$ be the higher order term in section 3.
Then we have
\bgn{lemma}
 $\Ob_k=\Ob_k(a_1,\cdots, a_{k-1}, b_1\cdots, b_{k-1})$
satisfies the following inequality,
$$
\|\Ob_k\|_{s-1}\leq C(\lam)M_k,
$$
where $C(\lam)$ depends on $\lam$ and we have
$$
\lim_{\lam\to 0}C(\lam)=0.
$$
\end{lemma}
\bgn{proof}
Since $\Ob_k$ determined by the terms of order greater than or equal to $2$,
$$\Ob_k=\sum_{l=2}^k \frac1{l!}(\ad_{z(t)}^l\, d)\psi.
$$%
We have
$$
\|\,[d, z(t)] \psi\,\|_{s-1}{\<} 2\|z(t)\psi\|_s.
$$%
Since $(\ad_{z(t)}^l\,d)=[\,\ad_{z(t)}^{l-1} \,\,d,\,z(t)\,]$, we find
\bgn{equation}
\|\(\ad_{z(t)}^l d\)_{[k]}\psi\|_{s-1}\leq 2(2C_s)^l(\|z(t)\|^l_s\, \|\psi\|_s)_{[k]}
\end{equation}
Hence it follows
\bgn{align}
\|\Ob_k\|_{s-1}=&\sum_{l=2}^k\frac1{l!}\big\|\(\ad_{z(t)}^l\, d\,\)_{[k]}\psi\big\|_{s-1} \\
\leq &\sum_{l=2}^k\frac1{l!}2(2C_s)^l\(\|z(t)\|_{s}^l\, \|\psi\|_{s}\)_{[k]}   \label{159}
\end{align}
Since the degree of $z(t)$ is greater than or equal to $1$, it follows from our assumption (\ref{156}) and $l\geq 2$ that we have
\bgn{equation}
\(\|z(t)\|_s^l\)_{[k]}\leq \(M(t)\)^l_{[k]}.\label{160}
\end{equation}
(Note that $\(\|z(t)\|_s^l\)_{[k]}$ consists of the term $\|z_i\|_s$, for $i<k$.) 
Substituting (\ref{160}) into (\ref{159}) and using lemma 4.2, we obtain
\bgn{align}
\|\Ob_k\|_{s-1}\leq &\sum_{l=2}^k\frac1{l!}2(2C_s)^l \(M(t)^l\)_{[k]}\|\psi\|_s ,\\
\leq &C_3\sum_{l=2}^k\frac1{l!}(2C_s)^l\lam^{l-1}M_k\\
\leq& C_3\lam^{-1}(e^{2C_s\lam}-1-2C_s\lam)M_k\\
=&C(\lam)M_k. \nonumber
\end{align}
where $C_3=2\|\psi\|_s$.
Then it follows the constant $C(\lam)$ satisfies
$$
\lim_{t\to 0} C(\lam)=0.
$$%
\end{proof}
\bgn{lemma}
$$\|b(t)\|_s\underset{k}{\<} K_2M(t).$$
\end{lemma}
\bgn{proof}
In section 3,
$b_k$ is defined as the solution of the equation,
\bgn{align}
\frac1{k!}db_k\psi+\frac1{k!}da_k\psi+\Ob_k=0
\end{align}
In fact $b_k\psi$ is given by
\bgn{align}
\frac1{k!}b_k\psi=-G_{K}d^*(\Ob_k)-G_{K}d^*
(\frac1{k!}a_k\psi)
\end{align}
Thus it follows from (\ref{149}) and the Schauder estimate
\bgn{align}
\|\frac1{k!}b_k\|_s\leq& C_2C_{K}\|\Ob_k\|_{s-1}+C_2C_{K}\|\frac1{k!}a_k\psi\|_s \label{166}
\end{align}
Applying lemma 4.3 and (\ref{142}) to (\ref{166}), we have
\bgn{align}
\|\frac1{k!}b_k\|_s
\leq &C_2C_{K}C(\lam)M_k+C_sC_2C_K K_1M_k\|\psi\|_s \nonumber\\
\leq&(C_4C(\lam)+C_5K_1)M_k \label{167}
\end{align}
where $C_4=C_2C_K$ and $C_5=C_sC_2\|\psi\|_s$.
Then from (\ref{150}) and (\ref{167}) if we take $K_2$ as
\bgn{align} \label{168}
K_2:=\max\{\, C_2C_1K_1,  (C_4C(\lam)+C_5K_1)\},
\end{align}
then we have  the inequality,
\bgn{equation}
\|b(t)\|_s\underset{k}{\<} K_2 M(t) \label{169}
\end{equation}
\end{proof}\\
Finally we estimate $z_k$.
It follows that
$$
(z(t))_{[k]}=\frac1{k!}z_k=\(e^{z(t)}-1-\sum_{p=2}^k \frac1{p!}z(t)^p\)_{[k]}.
$$
Hence we have
\bgn{equation}
\|\frac1{k!}z_k\|_s\leq \|(e^{z(t)}-1)_{[k]}\|_s+\sum_{p=2}^k
\frac1{p!}\|\(z(t)^p\)_{[k]}\|_s  \label{170}
\end{equation}
From our assumption and (\ref{169}),
$$
\|a(t)\|_s\<K_1M(t), \qquad \|b(t)\|_s\underset{k}{\<}K_2M(t).
$$
Then it follows from lemma \ref{lem: 4.1} and lemma \ref{lem: 4.2} that 

\bgn{align}\label{4.32}
\|e^{a(t)}-1\|_s&\< \frac1\lam (e^{K_1\lam} -1)M(t).
\end{align}
We also have
\bgn{align}\label{4.33}
\|e^{b(t)}-1\|_s&\underset{k}{\< }\frac 1\lam (e^{K_2\lam}-1)M(t)
\end{align}
Then we obtain
\bgn{lemma}
$$
\|z(t)\|_s\underset{k}{\<} M(t).
$$
\end{lemma}
\bgn{proof}
It follows from lemma \ref{lem: 4.2} that lemma \ref{142} that 
$$\|e^{a(t)}\|_s{\<}\frac 1{\lam}e^{K_1\lam} M(t).$$
Then substituting (\ref{4.32}) and (\ref{4.33}) into (\ref{4.34}), we have   

\bgn{align}
\|(e^{z(t)}-1)\|_s&\underset{k}{\<}
\|e^{a(t)}(e^{b(t)}-1)\|_s+\|e^{a(t)}-1\|_s \label{4.34}\\
&\underset{k}{\<}\frac 1\lam e^{K_1\lam} M(t)\frac 1\lam (e^{K_2\lam }-1)M(t) +\frac 1\lam (e^{K_1\lam} -1)M(t)\\
\end{align}
Applying lemma \ref{lem: 4.1} again, we have 
\bgn{align}
\|\(e^{z(t)}-1\)\|_s&\underset{k}{\<}\( e^{K_1\lam}\frac 1\lam(e^{K_2\lam}-1)  +\frac 1\lam (e^{K_1\lam}-1)  \)M(t)\\
&\underset{k}{\<}C(K_1, K_2)M(t)
\end{align}
where $C(K_1, K_2)$ is a constant which depends only on $K_1$ and $K_2$.
Since $(z(t))^p_{[k]}$ consists of terms $z_i$ for $i<k$, it follows from our assumption of the induction that
the second term of (\ref{170}) satisfies
\bgn{align}
\sum_{p=2}^k\frac1{p!}\|\(z(t)^p\)_{[k]}\|_s&\leq \sum_{p=2}^k\frac1{p!}\((C_sM(t))^p\)_{[k]}\\
&\leq \frac1{\lam}(e^{C_s\lam}-1-C_s\lam)M_k\\
&=C_1(\lam)M_k,
\end{align}
where $\lim_{\lam \to 0}C_1(\lam)=0.$
Thus if we take
$K_1, K_2$, $\lam$ which satisfy
\bgn{equation}   \label{180}
C(K_1, K_2)+C_1(\lam)\leq 1,
\end{equation}
 it follows from (\ref{170}) that
\bgn{align}
\frac1{k!}\|z_k\|_s\leq (C(K_1, K_2)+C_1(\lam))M_k\leq M_k
\end{align}
Thus
$\|z(t)\|_s\underset{k}{\<}M(t)$.
\end{proof}\\
If we take $a(t)$ sufficiently small,
we can take $K_1$, $K_2$ and $\lam$ with $K_1+K_2<1$ which satisfy
(\ref{168}) and (\ref{180}).
Hence by the induction, it turns out that $b(t)$ and $z(t)$ in section 3 are convergent series.

}

\large{
\section{Applications}

\subsection{\K\"ahler structures on K\"ahler manifolds}
Let $X$ be a compact K\"ahler manifold with the complex structure $J$ and the K\"ahler from $\ome$.
Then we have the \K\"ahler structure $(\J, e^{\sqrt{-1}\ome})$ with one pure spinor on $X$.
The deformations complex of \complex structures is given by the complex $(\w^\bullet \ol L, d_L)$.  
The complex
$(\w^\bullet \ol L, d_L)$ is isomorphic to the complex $(U^{-n+\bullet}\otimes K^{-1}_J, \pi_\bullet \circ d_{E_{0}})$, where $K^{-1}_J$ denotes the dual of the (usual) canonical line bundle of the complex manifold $(X, J)$.
In the case $(\J, e^{\sqrt{-1}\ome})$ on a K\"ahler manifold, we see that
$U^{-n+\bullet}$ is written in terms of the (usual) complex forms of type $(r,s)$,
\bgn{align}
&U^{-n}=\w^{n,0},\\
&U^{-n+1}=\w^{n,1}\oplus \w^{n-1,0},\\
&U^{-n+2}=\w^{n,2}\oplus\w^{n-1,1}\oplus\w^{n-2,0},\\
&U^{-n+3}=\w^{n,3}\oplus\w^{n-1,2}\oplus\w^{n-2,1}\oplus\w^{n-3,0}.
\end{align}
We take an open cover $\{V_\a\}$ of $X$ and $\Ome_\a$ as a nowhere vanishing holomorphic $n$-form on $V_\a$.
Then $E_{\a,0}=0$ and the operator $\pi_\bullet \circ d_{E_{\a,0}}$ is the (usual) $\ol\pa$ operator.
It implies that the space of infinitesimal deformations of \complex structures on $X$ is given by the direct sum of the $K^{-1}_J$-valued Dolbeault cohomology groups
\bgn{equation}\label{eq: 186}
H^{n,2}_{\ol\pa}(X, K^{-1}_J)\oplus H^{n-1,1}_{\ol\pa}(X, K^{-1}_J)\oplus H^{n-2,0}_{\ol\pa}(X, K^{-1}_J),
\end{equation}
where the space $H^{n-1,1}_{\ol\pa}(X, K^{-1}_J)\cong H^1(X,\Theta)$ is the space of
infinitesimal deformations of complex structures in Kodaira-Spencer theory.
The space $H^{n,2}_{\ol\pa}(X, K^{-1}_J)$ is given by the action of $B$-fields ($2$-forms) and
the space $H^{n-2,0}_{\ol\pa}(X, K^{-1}_J)$ is induced by the action of
holomorphic $2$-vector fields.

The space of the obstructions is given by
\bgn{align}
H^{n,3}_{\ol\pa}(X, K^{-1}_J)\oplus H^{n-1,2}_{\ol\pa}(X, K^{-1}_J)\oplus H^{n-2,1}_{\ol\pa}(X, K^{-1}_J)\oplus H^{n-3,0}_{\ol\pa}(X, K^{-1}_J).
\end{align}
Note that the description in equation (\ref{eq: 186}) is related to that in \cite{Gu1}.
Similarly we find that the first cohomology of the complex $(K^\bullet, d)$ is described as
\bgn{equation}
H^1(K^\bullet)\cong H^{1,1}_{\ol\pa}(X).
\end{equation}
Hence it follows from theorem 3.1 and 3.2 we obtain
\bgn{theorem}\label{theorem 5.1}
Let $X$ be a compact K\"ahler manifold with the \K\"ahler structure $(\J, e^{\sqrt{-1}\ome})$. If the obstruction space
$$\bigoplus_{i=0}^3 H^{n-i,3-i}_{\ol\pa}(X, K^{-1}_J)$$ vanishes, then we have the family of
\K\"ahler structures $\{ \J_t, \psi_{t,s}\}$ with $(\J_0, \psi_{\ss 0,0})=(\J, e^{\sqrt{-1}\ome})$ which is parametrized by $(t,s) \in \trian'\times W$,
where $\trian'$ is a small open set of $$\bigoplus_{i=0}^2 H^{n-i, 2-i}_{\ol\pa}(X, K^{-1}_J)$$
and $W$ denotes a small open set of  $H^{1,1}_{\ol\pa}(X)$ containing the origin.
\end{theorem}
There is no deformations of complex structures on the complex projective space $\C P^2$.
However there is a family of deformations of \complex structures on
$\C P^2$ which is parametrized by the space of holomorphic $2$-vector fields $H^0(\C P^2, \w^2\Theta)$.
Let $\{\, V_\a\, , \, \Ome_\a\}$ be a trivialization of the canonical line bundle $K$.
 Let $\b$ be a holomorphic $2$-vector field on $\C P^2$.
 Then it follows that the action of spin group on $\Ome_\a$
 $$
 e^{\b t}\w \Ome_\a
 $$
 induces deformations of \complex structure on $\C P^2$.
In fact, we take inhomogeneous coordinates $(z_1^\a, z_2^\a)$ on each $U_\a$ with $\Ome_\a=dz_1^\a\w dz_2^\a$,
and  $\b$ is written as
$$
\b=f\frac{\pa}{\pa z_1^\a}\w\frac{\pa}{\pa z_2^\a},
$$%
where $f$ is a cubic function.
Then
$$
e^{\b}\w \Ome_a=f+\Ome_\a.
$$%
Thus $e^{\b}\w \Ome_a$ is a non-degenerate pure spinor which
induces a \complex structure $\J_\b$. The type of \complex structure $\J$ is defined as the minimal degree of differential forms (non-degenerate pure spinors) which induces $\J$. Thus the type of $J_\b$ is 0 on the complement of the zero set of $\b$ and the type of $\J_\b$ is $2$ at
the zero set of $\b$.
Since we have $H^0(\C P^2,\w^2\Theta)\cong
H^0(\C P^2, {\mathcal O}(3))$, it follows from theorem of stability
that we have a family of \K\"ahler structures on $\C P^2$ parametrized by $H^0(\C P^2, {\mathcal O}(3))\oplus H^{1,1}_{\ol\pa}(X)$.
\subsection{\K\"ahler structures on Fano surfaces}
Our theorem can be applied to Fano surfaces.
Let $S_n$ be a blown up $\C P^2$ at $n$ points whose anti-canonical line bundle is ample ($n\leq 8$).
Then it follows from the Kodaira vanishing theorem that
the space of obstructions vanishes.
Thus
deformations
of \complex structures are parametrized by an open set of
$H^0(S_n, K^{-1})\oplus H^1(S_n, \Theta)$, whose dimensions are given by
$$\dim H^1(S_n, \Theta)=
\bgn{cases}& 2n-8,\quad (n =5,6,7,8),\\
& 0, \qquad\,\,\quad ( n=0,1,2,3,4)
\end{cases}
$$
$$\dim H^0(S_n, K^{-1})=
10-n
$$
It follows from theorem of stability we have the family of \K\"ahler structures on $S_n$ which is parametrized by an open set of the direct sum,
$$H^0(S_n, K^*)\oplus H^1(S_n, \Theta)\oplus H^{1,1}(S_n),
$$
where $H^{1,1}(S_n)$ denotes the Dolbeault cohomology of type $(1,1)$ which coincides with the cohomology $H^1(K^\bullet)$
(see section 4),
$$
\dim H^{1,1}(S_n)=1+n.
$$
\subsection{Poisson structures and \K\"ahler structures}
In general, we have an obstruction to deformations of \complex structures and the space of infinitesimal  deformations does not coincide with the space of actual deformations.
However theorem of stability can be applied as long as we have a one dimensional analytic family of deformations
of \complex structures.
Typical examples are constructed from holomorphic Poisson structures.
Let $X$ be a compact K\"ahler manifold with a holomorphic $2$-vector field $\b$.
If $\b$ satisfies that
\bgn{equation}
[\b, \b]_L=0,
\end{equation}
where the bracket denotes the Schouten bracket,
then $\b$ is called a holomorphic Poisson structure on $X$.
Since $\b$ is holomorphic, we find $d_L\b=0$.
Hence $\b$ also satisfies the Maurer-Cartan equation and the adjoint action of $e^{\b t}$ on $\J$
induces an analytic family of deformations of \complex structures.
We write it by $\J_{t\b}=\Ad_{e^{t\b}}\J$.
Hence we obtain from theorems 3.1 and 3.2
\bgn{theorem}\label{theorem 5.2}
Let $\b$ be a holomorphic Poisson structure on a compact K\"ahler manifold $X$.
Then we have a family of \K\"ahler structures $\{\J_{t\b}, \,\psi_t\}$.
\end{theorem}
The rank of $2$-vector $\b$ at $x$ is $r$ if $\b_x^r\neq 0$ and $\b^{r+1}_x=0$ for a point $x\in X$.  Then we denote it by rank $\b_x =r$.
 Since the type of \complex structure of $\J_\b$ is defined as the minimal degree of differential form $e^{\b}\cdot \Ome_\a$, where $\Ome_\a$ denotes a non-zero holomorphic $n$-form. Thus we have
 \bgn{equation}\label{189}
 \text{\rm type}(\J_\b)_x= n-2\, \text{\rm rank}\,\b_x.
 \end{equation}
This is concerned with the fact that the type $(\J_\b)_x$  can jump, depending on a choice of $x\in X$.
Let $X$ be a K\"ahler manifold with an action of an $l$ dimensional complex commutative Lie group $G$
 ($l\geq 2$).
We denote by $\{\xi_i\}_{i=1}^l$ a basis of the Lie algebra of $G$ which induces
 the corresponding holomorphic vector fields $\{V_i\}_{i=1}^l$ on $X$.
 We take $\b$ as a linear combination of $V_i\w V_j $'s,
 \bgn{equation}\label{190}
 \b=\sum_{i,j} \lam_{i,j} V_i\w V_j,
 \end{equation}
 where each $\lam_{i,j}$ denotes a constant.
 Since $[V_i, V_j]=0$, we have $[\b,\b]_L=0$.
 Then we have a family of \K\"ahler structure on $X$.
 The type of $\J_\b$ can change, according to the fixed points set of the action of $G$.
Hence we have
\bgn{theorem}\label{theorem 5.3}
 Let $X$ be a compact K\"ahler manifold of dimension $n$.
If we have an action of an $l$ dimensional complex commutative Lie group $G$ with a non-trivial
 $2$-vector $\b$ as in (\ref{190}), then we have
a family of deformations of nontrivial \K\"ahler structures on $X$.
\end{theorem}
 Since the type of $\J_\b$ is given by $n-2\,$rank$\,\b$ from (\ref{189}),
 it follows that generalized K\"ahler structures in theorem 5.3 are not obtained by
 the action of $B$-fields ($2$-forms) from usual K\"ahler structures.

Theorems \ref{theorem 5.1}, \ref{theorem 5.2} and \ref{theorem 5.3} imply that there are many examples of deformations of \K\"ahler structures on K\"ahler manifolds, such as
every toric K\"ahler manifolds and 
the Grassmannians. 
On a complex surface, any holomorphic section of anti-canonical bundle gives the Poisson structure.
There is a classification of holomorphic Poisson surfaces and 
we can count the dimensions of sections of anti-canonical bundles on a given holomorphic Poisson surfaces \cite{B.M}, 
\cite{Sa}.

}

\large{
\section{Appendix}
Let $\J$ be a \complex structure on a manifold $X$. Then we have the decomposition,
$$
(\TT)\otimes\C=L_\J\oplus \ol L_\J
$$%
We denote by $|a|$ the degree of $a\in \w^p\ol L_\J$, that is $p$. 
Then for $a \in \w^*\ol L_{\J}$, we define a graded bracket by 
$$
[d, a]_G= da-(-1)^{|a|}ad.
$$
We also define a bracket $[a,b]_L$ by 
$$
[a,b]_L= [d,a]_Sb-(-1)^{(|a|+1)|b|}b[d,a]_S.
$$
There is the following explicit description,
\bgn{proposition}
$[a,b]_L$ is an element of $\w^{|a|+|b|-1}\ol L_\J$ which is given in terms of the derived bracket,
\bgn{align}
[E_1\cdots E_n,\, F_1\cdots F_m]_S =& 
\sum_{i,j}(-1)^{i+j} E_1 \cdots\overset{\ch i}E_i \cdots E_n[E_i, F_j]_{d}F_1\cdots\overset{\ch j}F_j\cdots F_m 
\end{align}
for $E_i,\,F_j\in \ol L_J$, $i=1,\cdots, n,\,\,j=1,\cdots, m$.
\end{proposition} 
\bgn{proof}
The bracket $[a, b]_L$ is an operator acting on the differential forms $\w^*T^*$. Then it turns out that 
$$
[a,b]_S f\phi= f[a,b]_S\phi, \qquad \phi\in \w^*T^*.
$$
for a function $f$.
Thus $[a,b]_S$ is not a differential operator but an element of $\w^*\ol L_\J$.
Next we see that 
\bgn{align}
[E, F_1\cdots F_m]_S=& [\{ d, E\}, F_1\cdots F_m]_S\\
&=\sum_j (-1)^{j+1}[E, F_j]_{d}F_1\cdots\overset{\ch j}F_j\cdots F_m
\end{align}
Further for $a, b \in \w^*\ol L_J$ and $E\in \ol L_\J$, we have 
\bgn{equation}
[E\w a, \, b ]= a\w[E, b]_S- E[a, b]_S.
\end{equation}
Then by the induction, we have the result.
\end{proof}

}

\bgn{thebibliography}{99}
\bibitem{A.G.G}
V.~Apostolov, P.~Gauduchon, G.~Grantcharov,
{\it Bihermitian structures on complex surfaces},
Pro. London Math. Soc. {\bf 79} (1999), 414-428,
Corrigendum: {\bf 92}(2006), 200-202,  MR2192389, Zbl 1089.53503
\bibitem{B.K}
S.~Barannikov and M.~Kontsevich,
{\it Frobenius manifolds and formality of Lie algebras of
polyvector fields} (1998), International Math. Res. Notices, (4), 201-215, MR1609624, Zbl 0914.58004
\bibitem{B.G.C}
H.~Bursztyn, M~Gualtieri and G.~Cavalcanti,
{\it Reduction of Courant algebroids and generalized complex structures},
Adv. Math, {\bf 211}(2)(2007), 726-765, MR2323543, Zbl 1115.53056
\bibitem{B.M}
C.~ Bartocci, E.~ Macr$\grave{\text{\rm i}}$,
{\it Classification of Poisson surfaces},
math.AG/0402338, Commun. Contemp. Math. 7 (2005), no. 1, 89--95. MR2129789, Zbl 1071.14514
\bibitem{Ca}
G.~Cavalcanti,
{\it New aspects of $dd^c$-lemma},
Math.DG/0501406, D. Phil thesis
\bibitem{Ch}
C.C.~Chevalley,
{\it The algebraic theory of Spinors},
Columbia University Press, New York, 1954. viii+131 pp. MR0060497 
Zbl 0057.25901
\bibitem{Go1}
R.~Goto,
{\it Moduli spaces of topological calibrations,
Calabi-Yau, hyperK\"ahler, G$_2$ and Spin$(7)$ structures},
Internat. J. Math. 15 (2004), no. 3, 211--257. MR2060789, Zbl 1046.58002
\bibitem{Go2}
R.~Goto,
{\it On deformations of generalized Calabi-Yau, hyperK\"ahler, G$_2$ and Spin$(7)$ structures},
Math.DG/0512211
\bibitem{Go3}
R.~Goto,
{\it Poisson structures and \K\"ahler submanifolds},
 J. Math. Soc. Japan 61 (2009), no. 1, 107--132. MR2272873, Zbl 1160.53014
\bibitem{Gu1}
M.~Gualtieri,
{\it Generalized complex geometry}
,thesis, Math. DG/0401221, refined version, Math.DG/0703298
\bibitem{Gu2}
M.~Gualtieri,
{\it Hodge decomposition for generalized K\"ahler manifolds}, Math. DG/0409093,
Lecture at the String Theory and Geometry workshop, (2007) Oberwolfach
\bibitem{Hi1}
N.~Hitchin,
{\it Generalized Calabi-Yau manifolds},
Math. DG/0401221, Q. J. Math. 54 (2003), no. 3, 281--308. MR2013140, Zbl 1076.3201
\bibitem{Hi2}
N.Hitchin,
{\it  Instantons, Poisson structures and generalized K\"ahler geometry},
Comm. Math. Phys. 265 (2006), no. 1, 131--164. MR2217300, Zbl 1110.53056
\bibitem{Hi3}
N.Hitchin,
{\it Bihermitian metrics on Del Pezzo surfaces},
Math.DG/0608213,
J. Symplectic Geom. 5 (2007), no. 1, 1--8. MR2371181, Zbl pre05237677
\bibitem{Hu}
D.~Huybrechts,
{\it Generalized Calabi-Yau structures, K3 surfaces and B-fields}, math.AG/0306132,
Internat. J. Math. 16 (2005), no. 1, 13--36. MR2115675, Zbl 1120.14027
\bibitem{Ko}
K.~Kodaira
{\it Complex manifolds and deformations of complex structures},
Grundlehren der Mathematischen Wissenschaften [Fundamental Principles of Mathematical Sciences], 283. Springer-Verlag, New York, 1986. x+465 pp. ISBN: 0-387-96188-7 MR0815922, Zbl 0581.32012
Grundlehren der Mathematischen Wissenschaften, {\bf 283}, springer-Verlag,
(1986)
\bibitem{K.S.I,II}
K.~Kodaira and D.C.~Spencer,
{\it  On deformations of complex, analytic structures I,II},
Ann. of Math. (2) 67 (1958) 328--466. MR0112154, Zbl 0128.16901
\bibitem{K.S III}
K.~Kodaira and D.C.~Spencer,
{\it On deformations of complex analytic structure, III.
stability theorems for complex structures},
Ann. of Math. (2) 71 (1960) 43--76. MR0115189, Zbl 0128.16902 
\bibitem{Kos}
Y.~Kosmann-Schwarzbach
{\it Derived brackets}, 
Lett. Math. Phys. 69 (2004), 61--87. MR2104437, Zbl 1055.17016 
\bibitem{L.T}
Y.~Lin and S.~Tolman,
{\it Symmetries in generalized K\"ahler geometry},
Comm. Math. Phys. 268 (2006), no. 1, 199--222. MR2249799, Zbl 1120.53049
, Math. DG/0509069
\bibitem{L.W.P}
Z.-J.~Liu, A.~Weinstein and Ping. Xu, 
{\it Manin triples for Lie bialgebroids},
J. Differential Geom. 45 (1997), no. 3, 547--574. MR1472888, Zbl 0885.58030
\bibitem{Sa}
F.~Sakai,
{\it Anti-Kodaira dimension of ruled surfaces},
 Sci. Rep. Saitama Univ. Ser. A 10, no. 2, 1--7. (1982), MR0662405, Zbl 0496.14022
\bibitem{Se}
J. P. Serre,
{\it Lie Algebras and Lie groups},
Lecture Notes in Mathematics, 1500. Springer-Verlag, Berlin, 2006. viii+168 pp. ISBN: 978-3-540-55008-2; 3-540-55008-9 MR2179691, Zbl 0742.17008
\bibitem{Ti}
G.~Tian,
{\it Smoothness of the universal deformation space of compact Calabi-Yau manifolds and its Petersson-Weil metric}, Mathematical aspects of string theory (San Diego, Calif., 1986), 629--646, Adv. Ser. Math. Phys., 1, World Sci. Publishing, Singapore, 1987. MR0915841, Zbl 0696.53040
\end{thebibliography}
\bigskip
\bigskip
{Deaprtment of Mathematics \\ 
Graduate School of Science \\
Osaka University\\
Toyonaka, Osaka, 560\\
Japan
}
{goto@math.sci.osaka-u.ac.jp}

\end{document}